\documentclass[10pt]{article}

\usepackage[utf8]{inputenc}
\usepackage{amsmath}
\usepackage{amsfonts}
\usepackage{amssymb}
\usepackage{ntheorem}
\usepackage{pdfpages,textcomp}
\usepackage{graphicx}
\usepackage{makeidx,yfonts}
\usepackage{multicol}
\usepackage{enumerate,enumitem}
\usepackage{fontenc,mathtext}
\usepackage{array}
\usepackage{multirow}
\usepackage{float}
\usepackage{mathtools}
\restylefloat{table}
\usepackage[russian,english]{babel}

\newcommand{\HH}{\mathbb{H}}
\newcommand{\RR}{\mathbb{R}}

\theoremstyle{break}
\newtheorem{theorem}{Theorem}
\newtheorem{corollary}{Corollary} 
\newtheorem{proposition}{Proposition}    
\newtheorem{lemma}{Lemma}

\theorembodyfont{\upshape}
\theoremstyle{definition}
\newtheorem{definition}{Definition}   
\newtheorem{remark}{Remark}
\theoremstyle{nonumberplain}
\newtheorem{proof}{Proof}

\author{Matthieu Jacquemet\thanks{University of Applied Sciences and Arts Western Switzerland Valais and University of Fribourg. E-mail : jacquemet.matthieu@gmail.com} \hspace{0.5cm} Steven T Tschantz\thanks{Vanderbilt University, Nashville TN. E-mail : steven.tschantz@vanderbilt.edu}}
\title{All hyperbolic Coxeter $n$-cubes}
\date{\today}

\begin{document}
\maketitle

\begin{abstract}
Beside simplices, $n$-cubes form an important class of simple polyhedra. Unlike hyperbolic Coxeter simplices, hyperbolic Coxeter $n$-cubes are not classified. We show that there is no hyperbolic Coxeter $n$-cube for $n\geq~6$, and provide a full classification for $n\leq 5$. Our methods, which are essentially of combinatorial and algebraic nature, can be (and have been successfully) implemented in a symbolic computation software such as \textsf{Mathematica}$^\circledR$.
\end{abstract}

\section{Introduction}
Let $\HH^n$ be the $n$-dimensional real hyperbolic space. A \textit{hyperbolic Coxeter polyhedron} $\mathcal{P}\subset\HH^n$ is a finite-volume convex polyhedron whose angles are of the form $\frac{\pi}{k}$ for some $k\in\{2,...,\infty\}$. Identifying the facets of a hyperbolic Coxeter polyhedron by using the reflections in their supporting hyperplanes is a simple way to construct hyperbolic $n$-orbifolds and $n$-manifolds. In the known cases, such polyhedra are responsible for minimal volume hyperbolic orbifolds and manifolds~\cite{Kellerhals3}.\\

In contrast to the spherical and Euclidean cases, hyperbolic Coxeter polyhedra cannot exist any more if $n\geq 996$, and are far from being classified. In fact, comprehensive lists are available only if the number of facets of $\mathcal{P}$ equals $n+1$ or $n+2$. For example, hyperbolic Coxeter simplices exist only for $2\leq n\leq 9$. An overview of the current knowledge about the classification of hyperbolic Coxeter polyhedra (and related questions) is available on Anna Felikson's webpage~\cite{AnnaPage}, for example.\\
Hyperbolic $n$-cubes are simple polyhedra bounded by $2n$ facets in $\HH^n$, and, unlike simplices, they have no simplex facet. Hyperbolic space forms can be constructed by identifying isometric facets of hyperbolic $n$-cubes in a suitable way. Notably, hyperbolic manifolds which can be decomposed into regular ideal cubes (so-called \textit{cubical manifolds}) are not necessarily decomposable into regular ideal tetrahedra.\\

In a previous paper \cite{Jacquemet3}, the first author showed that there are no hyperbolic Coxeter $n$-cubes in $\HH^n$ for $n\geq 10$, and a combinatorial/algebraic approach has been used in order to classify the \textit{ideal} hyperbolic Coxeter $n$-cubes, which exist only for $n\leq 3$. As a consequence of the relatively low dimensional upper bound and the explicit approach leading to the classification in the ideal case, it was finally suggested that a computer assisted approach could successfully handle the large amount of combinatorial possibilities in the general setting, and therefore proceed with the classification.\\
In this paper, we follow this strategy and eventually provide a full classification of \textit{all} hyperbolic Coxeter $n$-cubes. To this end, we refine and expand the combinatorial and algebraic methods from~\cite{Jacquemet3}, and we implement them with help of \textsf{Mathematica}.\\
First, we enumerate all \textit{potential hyperbolic Coxeter $n$-cube matrices}, that is, all Coxeter matrices satisfying a certain set of necessary conditions arising from the combinatorial setting. This already allows us to show that there is no hyperbolic Coxeter $n$-cube in dimensions $n\geq 7$. Then, we consider the polynomial (in)equalities which have to be satisfied by entries of these potential matrices in order to extract from that list the matrices which lead to all the Gram matrices of hyperbolic Coxeter $n$-cubes. A key tool in this context is the extraction of \textit{Gröbner bases} of sets of polynomials.\\

More precisely, we show the following result.
\begin{theorem}[Classification of the hyperbolic Coxeter $n$-cubes]\label{thm classification}
The hyperbolic Coxeter $n$-cubes are repartitioned as follows.
\begin{enumerate}
\item The hyperbolic Coxeter 2-cubes and 3-cubes are distributed according to Tables~\ref{table 2c} and \ref{table 3c}, respectively. An explicit list of all hyperbolic 3-cubes families can be found in the Appendix, p.22.
\item There are exactly 15 hyperbolic Coxeter 4-cubes: 12 compact ones, whose graphs are given in Figure~\ref{fig 4cubes}, and 3 noncompact ones, whose graphs are given in Figure~\ref{fig 4cubesnc}.
\item There are exactly 4 hyperbolic Coxeter $5$-cubes: a unique compact one, whose graph is given in Figure~\ref{fig 5cubec}, and 3 noncompact ones, whose graphs are given in Figure~\ref{fig 5cubesnc}.
\item There are no hyperbolic Coxeter $n$-cubes for $n\geq 6$.
\end{enumerate}
\end{theorem}

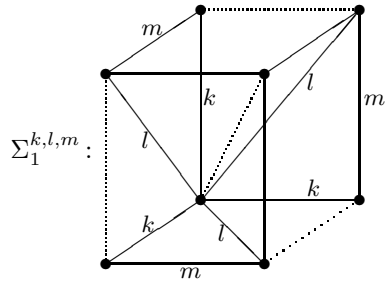
\begin{table}[H]
 \begin{center}\renewcommand{\arraystretch}{1.25}
 \begin{tabular}{| m{0.35\linewidth} || m{0.05\linewidth} | m{0.05\linewidth} |m{0.05\linewidth} |m{0.05\linewidth} |m{0.05\linewidth} |}
    \hline
    \centering \text{Number of free parameters}
    & \centering 0
    & \centering 1
    & \centering 2
    & \centering 3
    & \centering 4
	\tabularnewline \hline\hline
	\centering \text{Number of families}
    & \centering 119
    & \centering 75
    & \centering 30
    & \centering 5
    & \centering 1
     \tabularnewline\hline
 \end{tabular}
 \caption{Repartition of the hyperbolic Coxeter 2-cubes}\label{table 2c}
\end{center}
\end{table}
 
\begin{table}[H]
 \begin{center}\renewcommand{\arraystretch}{1.25}
 \begin{tabular}{| m{0.35\linewidth} || m{0.08\linewidth} | m{0.08\linewidth} |m{0.08\linewidth} |m{0.08\linewidth} |}
    \hline
    \centering \text{Number of free parameters}
    & \centering 0
    & \centering 1
    & \centering 2
    & \centering 3
	\tabularnewline \hline\hline
	\centering \text{Compact}
    & \centering 2016
    & \centering 315
    & \centering 17
    & \centering 1
    \tabularnewline \hline
    \centering \text{Noncompact}
    & \centering 6250
    & \centering 533
    & \centering 10
    & \centering 0
     \tabularnewline\hline
 \end{tabular}
 \caption{Repartition of the hyperbolic Coxeter $3$-cubes}\label{table 3c}
\end{center}
\end{table}

\begin{figure}[h!]
\begin{center}
\setlength{\unitlength}{0.8pt}
\begin{picture}(350,340)
	\put(120,210){\circle*{5}}
	\put(165,240){\circle*{5}}
	\put(120,300){\circle*{5}}
	\put(165,330){\circle*{5}}
	\put(195,210){\circle*{5}}
	\put(240,240){\circle*{5}}	
	\put(195,300){\circle*{5}}
	\put(240,330){\circle*{5}}	
	\multiput(120,210)(0,2){45}{\line(0,1){0.5}}
	\multiput(165,330)(2,0){38}{\line(0,1){0.5}}
	\multiput(195,210)(3,2){15}{\line(0,1){0.5}}
	\multiput(165,240)(1,2){30}{\line(0,1){0.5}}
	\put(120,210){\line(3,2){45}}
	\put(120,210){\line(1,0){75}}
	\put(165,240){\line(1,-1){30}}
	\put(165,240){\line(1,0){75}}
	\put(165,240){\line(-3,4){45}}
	\put(165,240){\line(0,1){90}}
	\put(165,240){\line(5,6){75}}
	\put(195,300){\line(0,-1){90}}
	\put(195,300){\line(3,2){45}}
	\put(195,300){\line(-1,0){75}}
	\put(120,300){\line(3,2){45}}	
	\put(240,330){\line(0,-1){90}}
	\put(75,260){$\Sigma^{k,l,m}_{1}$\,:}
	\put(155,203){\small $m$}
	\put(242,285){\small $m$}
	\put(137,318){\small $m$}
	\put(137,225){\small $k$}
	\put(215,242){\small $k$}
	\put(166,285){\small $k$}
	\put(215,292){\small $l$}
	\put(137,265){\small $l$}
	\put(173,221){\small $l$}
	\put(60,180){\small $(k,l)\in\{(2,2),(2,3),(3,2)\}$}
	\put(220,180){\small $m=3,4,5$}
	
	\put(30,10){\circle*{5}}
	\put(30,130){\circle*{5}}
	\put(150,10){\circle*{5}}
	\put(150,130){\circle*{5}}
	\put(60,40){\circle*{5}}
	\put(120,40){\circle*{5}}	
	\put(90,70){\circle*{5}}
	\put(90,100){\circle*{5}}	
	\multiput(30,10)(2,2){15}{\line(0,1){0.5}}
	\multiput(150,10)(-2,2){15}{\line(0,1){0.5}}
	\multiput(30,130)(2,0){60}{\line(2,0){0.5}}
	\multiput(90,70)(0,2){15}{\line(0,2){0.5}}
	\put(30,10){\line(0,1){120}}
	\put(30,10){\line(1,0){120}}
	\put(150,10){\line(0,1){120}}
	\put(90,70){\line(-1,-1){30}}
	\put(90,70){\line(1,-1){30}}
	\put(90,70){\line(1,1){60}}
	\put(90,70){\line(-1,1){60}}
	\put(90,100){\line(-1,-2){30}}
	\put(90,100){\line(1,-2){30}}
	\put(90,100){\line(2,1){60}}
	\put(90,100){\line(-2,1){60}}
	\put(0,70){$\Sigma^{m}_{2}$\,:}
	\put(35,70){\small $4$}
	\put(140,70){\small $4$}
	\put(65,113){\small $m$}
	\put(65,68){\small $m$}
	\put(107,113){\small $m$}
	\put(107,68){\small $m$}
	\put(76,-10){\small $m=2,3$}

	\put(220,10){\circle*{5}}
	\put(250,40){\circle*{5}}
	\put(220,100){\circle*{5}}
	\put(250,130){\circle*{5}}
	\put(310,10){\circle*{5}}
	\put(340,40){\circle*{5}}	
	\put(310,100){\circle*{5}}
	\put(340,130){\circle*{5}}	
	\multiput(220,10)(0,2){45}{\line(0,1){0.5}}
	\multiput(310,10)(0,2){45}{\line(0,1){0.5}}
	\multiput(250,40)(2,0){45}{\line(2,0){0.5}}
	\multiput(250,130)(2,0){45}{\line(2,0){0.5}}
	\put(220,10){\line(1,1){30}}
	\put(220,100){\line(1,1){30}}
	\put(310,10){\line(1,1){30}}
	\put(310,100){\line(1,1){30}}
	\put(220,10){\line(1,0){90}}
	\put(220,100){\line(1,0){90}}
	\put(250,40){\line(0,1){90}}
	\put(340,40){\line(0,1){90}}
	\put(195,70){$\Sigma_{3}$\,:}
	\put(230,26){\small $4$}
	\put(230,116){\small $4$}
	\put(320,26){\small $4$}
	\put(320,116){\small $4$}
\end{picture}
\end{center}
\caption{The graphs of the compact hyperbolic Coxeter $4$-cubes (the weights of the dotted edges are given in Section~\ref{subsec dim4})}\label{fig 4cubes}
\end{figure}

\begin{figure}[h!]
\begin{center}
\setlength{\unitlength}{1.1pt}
\begin{picture}(235,110)
	\put(30,10){\circle*{3}}
	\put(60,10){\circle*{3}}
	\put(90,10){\circle*{3}}
	\put(30,40){\circle*{3}}
	\put(60,40){\circle*{3}}
	\put(90,40){\circle*{3}}
	\put(45,70){\circle*{3}}
	\put(75,70){\circle*{3}}
	\multiput(30,10)(0,2){15}{\line(0,1){0.5}}
	\multiput(60,10)(0,2){15}{\line(0,1){0.5}}
	\multiput(90,10)(0,2){15}{\line(0,1){0.5}}
	\multiput(45,70)(2,0){15}{\line(1,0){0.5}}
	\put(30,10){\line(1,1){30}}
	\put(30,10){\line(2,1){60}}
	\put(30,40){\line(1,-1){30}}
	\put(30,40){\line(2,-1){60}}
	\put(60,10){\line(1,1){30}}
	\put(60,40){\line(1,-1){30}}
	\put(30,10){\line(1,1){30}}
	\put(30,40){\line(1,2){15}}
	\put(30,40){\line(3,2){45}}
	\put(60,40){\line(-1,2){15}}
	\put(60,40){\line(1,2){15}}
	\put(90,40){\line(-3,2){45}}
	\put(90,40){\line(-1,2){15}}
	\put(0,40){$\Sigma^{m}_{4}$\,:}
	\put(43,45){\small $m$}
	\put(65,46){\small $m$}
	\put(82,56){\small $m$}
	\put(48,-5){\small $m=2,3$}
	
	\put(145,25){\circle*{3}}
	\put(175,25){\circle*{3}}
	\put(205,25){\circle*{3}}
	\put(235,25){\circle*{3}}
	\put(145,55){\circle*{3}}
	\put(175,55){\circle*{3}}
	\put(205,55){\circle*{3}}
	\put(235,55){\circle*{3}}

	\multiput(145,25)(0,2){15}{\line(0,1){0.5}}
	\multiput(175,25)(0,2){15}{\line(0,1){0.5}}
	\multiput(205,25)(0,2){15}{\line(0,1){0.5}}
	\multiput(235,25)(0,2){15}{\line(0,1){0.5}}
	\put(145,25){\line(1,1){30}}
	\put(145,25){\line(2,1){60}}
	\put(145,25){\line(3,1){90}}
	\put(145,55){\line(1,-1){30}}
	\put(145,55){\line(2,-1){60}}
	\put(145,55){\line(3,-1){90}}
	\put(175,25){\line(1,1){30}}
	\put(175,25){\line(2,1){60}}
	\put(175,55){\line(1,-1){30}}
	\put(175,55){\line(2,-1){60}}
	\put(205,25){\line(1,1){30}}
	\put(205,55){\line(1,-1){30}}

	\put(120,40){$\Sigma_5$\,:}
\end{picture}
\end{center}
\caption{The graphs of the noncompact hyperbolic Coxeter $4$-cubes (the weights of the dotted edges are given in Section~\ref{subsec dim4})}\label{fig 4cubesnc}
\end{figure}

\begin{figure}[h!]
\begin{center}
\setlength{\unitlength}{0.8pt}
\begin{picture}(150,135)
	\put(30,0){\circle*{4}}
	\put(110,0){\circle*{4}}
	\put(30,80){\circle*{4}}
	\put(70,40){\circle*{4}}
	\put(70,120){\circle*{4}}
	\put(110,80){\circle*{4}}
	\put(150,40){\circle*{4}}
	\put(150,120){\circle*{4}}
	\put(60,70){\circle*{4}}
	\put(120,50){\circle*{4}}
	\multiput(30,0)(2,0){40}{\line(1,0){0.5}}
	\multiput(70,40)(0,2){40}{\line(0,1){0.5}}
	\multiput(110,80)(2,2){20}{\line(0,1){0.5}}
	\multiput(150,40)(-3,1){10}{\line(0,1){0.5}}
	\multiput(30,80)(3,-1){10}{\line(0,1){0.5}}
	\put(30,0){\line(0,1){80}}
	\put(30,0){\line(1,1){40}}
	\put(70,120){\line(-1,-1){40}}
	\put(70,120){\line(1,0){80}}
	\put(110,80){\line(-1,0){80}}
	\put(110,80){\line(0,-1){80}}
	\put(150,40){\line(-1,0){80}}
	\put(150,40){\line(0,1){80}}
	\put(150,40){\line(-1,-1){40}}
	\put(60,70){\line(3,-1){60}}
	\put(0, 50){$\Gamma_1$\,:}
	
\end{picture}
\end{center}
\caption{The graph of the unique compact hyperbolic Coxeter $5$-cube (the weights of the dotted edges are given in Section~\ref{subsec dim5})}\label{fig 5cubec}
\end{figure}

\begin{figure}[h!]
\begin{center}
\setlength{\unitlength}{0.8pt}
\begin{picture}(320,170)
	\put(30,15){\circle*{4}}
	\put(110,15){\circle*{4}}
	\put(30,95){\circle*{4}}
	\put(70,55){\circle*{4}}
	\put(70,135){\circle*{4}}
	\put(110,95){\circle*{4}}
	\put(150,55){\circle*{4}}
	\put(150,135){\circle*{4}}
	\put(90,75){\circle*{4}}
	\put(130,75){\circle*{4}}
	\multiput(30,15)(2,0){40}{\line(1,0){0.5}}
	\multiput(70,55)(0,2){40}{\line(0,1){0.5}}
	\multiput(110,95)(2,2){20}{\line(0,1){0.5}}
	\multiput(150,55)(-2,2){10}{\line(0,1){0.5}}
	\multiput(30,95)(3,-1){20}{\line(0,1){0.5}}
	\put(30,15){\line(0,1){80}}
	\put(30,15){\line(1,1){40}}
	\put(70,135){\line(-1,-1){40}}
	\put(70,135){\line(1,0){80}}
	\put(110,95){\line(-1,0){80}}
	\put(110,95){\line(0,-1){80}}
	\put(150,55){\line(-1,0){80}}
	\put(150,55){\line(0,1){80}}
	\put(150,55){\line(-1,-1){40}}
	\put(130,75){\line(-1,0){40}}
	\put(130,75){\line(1,3){20}}
	\put(130,75){\line(-1,-3){20}}
	\put(130,75){\line(-3,-1){60}}
	\put(-5,70){$\Gamma^{m}_{2}$\,:}
	\put(50,41){\small $4$}
	\put(105,137){\small $4$}
	\put(104,41){\small $4$}
	\put(120,41){\small $m$}
	\put(129,97){\small $m$}
	\put(98,76){\small $m$}
	\put(82,62){\small $m$}
	\put(55,0){\small $m=2,3$}
	\put(24,17){\tiny 2}
	\put(104,17){\tiny 1}
	\put(24,97){\tiny 7}
	\put(104,97){\tiny 6}
	\put(64,137){\tiny 4}
	\put(144,137){\tiny 5}
	\put(84,78){\tiny 8}
	\put(124,77){\tiny 9}
	\put(152,57){\tiny 10}
	\put(64,57){\tiny 3}
	
	\put(220,15){\circle*{4}}
	\put(300,15){\circle*{4}}
	\put(200,95){\circle*{4}}
	\put(320,95){\circle*{4}}
	\put(240,35){\circle*{4}}
	\put(280,35){\circle*{4}}
	\put(220,85){\circle*{4}}
	\put(300,85){\circle*{4}}
	\put(260,115){\circle*{4}}
	\put(260,135){\circle*{4}}
	\multiput(200,95)(2,-1){10}{\line(0,1){0.5}}
	\multiput(220,15)(2,2){10}{\line(0,1){0.5}}
	\multiput(300,15)(-2,2){10}{\line(0,1){0.5}}
	\multiput(320,95)(-2,-1){10}{\line(0,1){0.5}}
	\multiput(260,135)(0,-2){10}{\line(0,1){0.5}}
	\put(200,95){\line(3,2){60}}
	\put(200,95){\line(1,-4){20}}
	\put(300,15){\line(1,4){20}}
	\put(300,15){\line(-1,0){80}}
	\put(320,95){\line(-3,2){60}}
	\put(240,35){\line(1,4){20}}
	\put(280,35){\line(-1,4){20}}
	\put(220,85){\line(1,0){80}}
	\put(220,85){\line(6,-5){60}}
	\put(300,85){\line(-6,-5){60}}
	\put(175,70){$\Gamma_{3}$\,:}	
\end{picture}
\end{center}
\caption{The graphs of the noncompact hyperbolic Coxeter $5$-cubes (the weights of the dotted edges are given in Section~\ref{subsec dim5})}\label{fig 5cubesnc}
\end{figure}

\pagebreak
\section{Hyperbolic Coxeter $n$-cubes}\label{sec cubes}

We denote by $\HH^n$ the hyperbolic space of dimension $n$, and by $\partial\HH^n$ its boundary. We set $\overline{\HH^n}=\HH^n\cup\partial\HH^n$.

\subsection{Hyperbolic Coxeter polyhedra}\label{subsec polyhedra}
Let $\mathbb{X}^n\in\{\mathbb{S}^n,\mathbb{E}^n,\HH^n\}$ be one of the three standard geometric spaces of constant curvature, realized in a suitable linear real space. A \textit{Coxeter ($n$-)polyhedron} is a convex, finite-volume $n$-polyhedron $\mathcal{P}\subset\mathbb{X}^n$ whose dihedral angles are of the form $\frac{\pi}{k}$, for $k\in\{2,...,\infty\}$. Standard references about Coxeter polyhedra and their properties are \cite{Ratcliffe,Vinberg2}.\\

In the sequel, we assume that $\mathcal{P}\subset\HH^n$ is a hyperbolic Coxeter polyhedron. Then, it is bounded by finitely many hyperplanes, say $H_1,...,H_N$, $N\geq n+1$, each $H_i$ being associated to a normal unit vector $u_i$. For $i=1,...,N$, the \textit{facet (or $(n-1)$-face)} $F_i$ of $\mathcal{P}$ is the intersection $F_i=\mathcal{P}\cap H_i$. For $0\leq k\leq n-2$, a \textit{$k$-face of $\mathcal{P}$} is a facet of a $(k+1)$-face of $\mathcal{P}$. A \textit{vertex} is a $0$-face, and an \textit{edge} is a $1$-face of $\mathcal{P}$.\\
If a vertex $v$ of $\mathcal{P}$ lies on $\partial\HH^n$ we call $v$ an \textit{ideal vertex} of $\mathcal{P}$. Then, $\mathcal{P}$ is a \textit{noncompact} polyhedron. If all vertices of $\mathcal{P}$ lie on $\partial\HH^n$, then $\mathcal{P}$ is said to be \textit{ideal}.\\

Let $v\in\HH^n$ be an ordinary vertex of $\mathcal{P}$. The \textit{vertex figure, or link,} $L(v)$ of $v$ is the intersection
\[L(v)=\mathcal{P}\cap\mathcal{S}_\rho(v),\]
where $\mathcal{S}_\rho(v)$ is a sphere with center $v$ and radius $\rho>0$ not containing any other vertex of $\mathcal{P}$ and not intersecting any facet of $\mathcal{P}$ not incident to $v$. In particular, $L(v)$ is a spherical Coxeter $(n-1)$-polyhedron.\\
Similarly, the link $L(v)$ of an ideal vertex $v\in\partial\HH^n$ is defined as the intersection of $\mathcal{P}$ with a sufficiently small horosphere centred at $v$. In particular, $L(v)$ is a Euclidean Coxeter $(n-1)$-polyhedron \cite[I, Chapter 6.2]{Vinberg2}.\\

The \textit{Gram matrix} of $\mathcal{P}$ is the matrix $G=G(\mathcal{P})=(g_{ij})_{1\leq i,j\leq N}$ given by
\[g_{ij}=\left\{\begin{array}{cl} 
1&,\,\text{if } j=i,\\ 
-\cos\frac{\pi}{k_{ij}} &,\, \text{if }H_i\text{ and }H_j\text{ intersect in }\HH^n \text{ with angle }\frac{\pi}{k_{ij}},\\ 
-1&,\, \text{if } H_i \text{ and } H_j\text{ are parallel},\\
-\cosh\,l_{ij} &,\, \text{if }d(H_i,H_j)=l_{ij}>0.
\end{array}\right.\]
The matrix $G$ is real, symmetric, and of signature $(n,1)$ \cite[Chapter 6.2]{Vinberg2}. If $\mathcal{P}$ is realized in the vector space model $\mathcal{H}^n$ for $\mathbb{H}^n$ and has normal vectors $u_1,...,u_N$, then $g_{ij}= \langle u_i,u_j\rangle$ for all $i,j\in\{1,...,N\}$, where $\langle .,.\rangle$ is the standard bilinear form of signature $(n,1)$ on $\mathcal{H}^n$.\\

A Coxeter polyhedron $\mathcal{P}\subset\HH^n$ is often described by its \textit{Coxeter graph} $\Gamma=\Gamma(\mathcal{P})$ as follows. A node $i$ in $\Gamma$ represents the bounding hyperplane $H_i$ of $\mathcal{P}$. Two nodes $i$ and $j$ are joined by an edge with weight $2\leq k_{ij}\leq\infty$ if $H_i$ and $H_j$ intersect in $\overline{\HH^n}$ with angle $\frac{\pi}{k_{ij}}$. If the hyperplanes $H_i$ and $H_j$ admit a common perpendicular of length $l_{ij}>0$ in $\mathbb{H}^n$, the nodes $i$ and $j$ are joined by a dotted edge, sometimes labelled $\cosh\,l_{ij}$. In practice and in the following discussion, an edge of weight $2$ is omitted, and an edge of weight $3$ is written without its weight. The \textit{rank} of $\Gamma$ denotes the number of its nodes.\\
The Coxeter graphs of indecomposable spherical (resp. Euclidean) Coxeter polyhedra are well known. These polyhedra exist in any dimension and are completely classified. The corresponding graphs are called \textit{elliptic} (resp. \textit{parabolic}). They can be found in \cite[pp. 202-203]{Vinberg2}, for example.\\

The \textit{Coxeter matrix} of $\mathcal{P}$ is the symmetric matrix $M=(m_{ij})_{1\leq i,j\leq N}$ with entries in $\mathbb{N}\cup\{\infty\}$ such that
\[m_{ij}=\left\{\begin{array}{cl} 
1&,\,\text{if } j=i,\\ 
k_{ij} &,\,\text{if }H_i\text{ and }H_j\text{ intersect in }\HH^n\text{ with angle }\frac{\pi}{k_{ij}},\\
\infty &,\, \text{otherwise}.
\end{array}\right.\]
Notice that the Coxeter matrix does not distinguish between parallel and ultra-parallel pairs of facets of $\mathcal{P}$.

\begin{remark}
In the sequel, we shall refer to \textit{the Coxeter matrix $M$ of a graph $\Gamma$} as the Coxeter matrix $M$ of the Coxeter polyhedron $\mathcal{P}$ such that $\Gamma=\Gamma(\mathcal{P})$.\\
Similarly, we say that a Coxeter graph or a Coxeter matrix is \textit{of (non)compact type} if the associated Coxeter polyhedron is (non)compact.
\end{remark}

\subsection{Hyperbolic Coxeter $n$-cubes}\label{subsec gencube}

A hyperbolic \textit{$n$-cube}, $n\geq 2$, is a polyhedron $\mathcal{C}\subset \HH^n$ whose closure $\overline{\mathcal{C}}\subset\overline{\HH^n}$ is combinatorially equivalent to the standard cube $[0,1]^n\subset\RR^n$. In particular, an $n$-cube has $2^n$ vertices, and is bounded by $n$ pairs of mutually disjoint hyperplanes. It is \textit{cubical}, i.e. its $k$-faces are $k$-cubes, $2\leq k\leq n-1$. Moreover, the number $f_k(\mathcal{C})$ of $k$-faces of $\mathcal{C}$ is given by (see \cite[Chapter 4.4]{Grun}, for example)
\[f_k(\mathcal{C})=2^{n-k}\begin{pmatrix}n\\k\end{pmatrix},\quad 0\leq k\leq n.\]

For $n\geq 2$, let $\mathcal{C}\subset\HH^n$ be an $n$-cube bounded by hyperplanes $H_1,...,H_{2n}$ such that the hyperplane $H_{2i-1}$ intersects all hyperplanes except $H_{2i}$ for $i=1,...,n$. This is going to be our preferred facet labelling in the sequel. The set $\mathfrak{H}=\{H_1,...,H_{2n}\}$ can be partitioned into two families of $n$ concurrent hyperplanes in $2^{n-1}$ different ways. Let $\mathfrak{H}=\mathfrak{H}_1\sqcup\mathfrak{H_2}$ be such a partition. Then, for $i=1,2$, the hyperplanes in $\mathfrak{H}_i$ form a simplicial cone in $\HH^n$ based at a vertex of $\mathcal{C}$, say $\mathfrak{v}_i$. The vertices $\mathfrak{v}_1$ and $\mathfrak{v}_2$ lie on a (spatial) diagonal of $\mathcal{C}$ (they are opposite in~$\mathcal{C}$).\\
Moreover, in the Coxeter case, for any vertex $p_i\in\mathcal{C}$, the graph of the vertex figure $L(p_i)$ is the subgraph of $\Gamma(\mathcal{C})$ of rank $n$ spanned by the nodes representing the hyperplanes in $\mathfrak{H}$ which contain $p_i$.

\section{Potential hyperbolic Coxeter $n$-cube matrices}\label{sec potential}
The first step of the classification is the enumeration of all matrices that are candidates for being Coxeter matrices of hyperbolic Coxeter $n$-cubes. We start with the concepts of partial matrices and potential matrices, which are going to allow a systematic implementation of the procedure.

\subsection{Partial matrices}\label{subsec partial}

\begin{definition}
Let $\Omega=\{n\,|\,n\geq 2\}\cup\{\infty\}$ and let $\bigstar$ be a symbol representing an undetermined real value. A \textit{partial matrix of size $m\geq 1$} is a symmetric $m\times m$ matrix $M$ whose diagonal entries are $1$, and whose non-diagonal entries belong to $\Omega\cup\{\bigstar\}$.
\end{definition}

In particular, a partial matrix $M$ with no $\bigstar$ entries is totally determined. Otherwise, there will be various ways to fill in undetermined values (the entries denoted by $\bigstar$ need not be equal).

\begin{definition}
Let $M=(m_{ij})_{1\leq i,j\leq m}$ and $M'=(m'_{ij})_{1\leq i,j\leq m}$ be partial matrices of size $m$. We say that $M'$ \textit{refines} $M$ if in every position $(i,j)$ where $m_{ij}\neq\bigstar$, one has $m'_{ij}=m_{ij}$, that is, $M'$ has all the same (determined) values as in $M$, and possibly more.
\end{definition}

Fix a total order on the entries in matrices of size $m$. If $M=(m_{ij})_{1\leq i,j\leq m}$ and $M'=(m'_{ij})_{1\leq i,j\leq m}$ are partial matrices, we write $M<M'$ if the first position $(i,j)$ in which the matrices differ has $m_{ij}<m'_{ij}$.\\
We write $M\leq M'$ if $M<M'$ or $M=M'$. Note that $\leq$ is a total order (a lexicographic ordering reading matrix entries in a specified way).\\
We define a partial ordering of partial matrices by $M\prec M'$ if some $m_{ij}<m'_{ij}$ are not $\bigstar$ and at all earlier positions $(i,j)$, one has $m_{ij}=m'_{ij}\neq \bigstar$, i.e. the matrices are determined (and coincide) up to and including the first position where they differ, where one has $m_{ij}<m'_{ij}$. Notice that if $M$ and $M'$ are totally determined, then $M\prec M'$ if and only if $M<M'$.\\
The following result is immediate.

\begin{proposition}
Let $M,M',N$ and $N'$ be partial matrices of size $m$. If $M'$ refines $M$ and $N'$ refines $N$, then $M\prec N$ implies $M'\prec N'$.
\end{proposition}

\begin{definition}
Let $M$ be an arbitrary $m\times m$ matrix, and $s=(s_1,s_2,...,s_k)$, $1\leq k\leq m$, be a sequence of distinct indices from $\{1,...,m\}$. Let $M^{s}$ be the $k\times k$ submatrix of $M$ with $(i,j)$-entry $m_{s_i,s_j}$. If $\sigma$ is a permutation of $\{1,...,l\}$, $1\leq l\leq m$, we write $M^\sigma$ for $M^{(\sigma(1),\sigma(2),...,\sigma(l))}$.
\end{definition}

\subsection{Potential matrices}\label{subsec potential}

If $M$ is the Coxeter matrix of a Coxeter $n$-polyhedron $\mathcal{P}$ with $m\geq n+1$ facets with respective normal vectors $u_1,...,u_m$, and if $\sigma$ is a permutation of $\{1,...,m\}$, then the polyhedron with normal vectors $u_{\sigma(1)},...,u_{\sigma(m)}$ is an isometric copy of $\mathcal{P}$, with facets relabeled and Coxeter matrix $M^\sigma$. Thus we will consider $M$ and $M^\sigma$ equivalent, for the purpose of understanding Coxeter polyhedra without preferred labelling. It suffices then to consider the representative $M$ such that $M\leq M^\sigma$ for every permutation $\sigma$ of $\{1,...,m\}$.\\

As mentioned in Section \ref{subsec polyhedra}, if an abstract polyhedron is to be realized as a hyperbolic Coxeter polyhedron with Coxeter matrix $M$, then for each subset of facets meeting in a vertex, with $s$ a sequence of indices of these facets, then $M^s$ is the Coxeter matrix for a finite Coxeter group (if the vertex is an ordinary vertex), or for an affine Coxeter group (if the vertex is an ideal vertex).\\

Recall that for any hyperbolic Coxeter polyhedron $\mathcal{P}\subset\HH^n$, any two hyperplanes bounding $\mathcal{P}$ intersect in $\HH^n$ if and only if their intersection is an $(n-2)$-face of $\mathcal{P}$ and the ridge of a dihedral angle of $\mathcal{P}$ \cite{Vinberg2}. 
\begin{definition}
We say that a Coxeter matrix $M$ is a \textit{potential matrix for a (polyhedron of a) given combinatorial type} if there are entries $\infty$ in positions in $M$ corresponding to non-adjacent facets, and for every sequence $s$ of indices of facets meeting at a vertex, the submatrix $M^s$ is the Coxeter matrix of a finite or affine Coxeter group. In this case, we say that the associated abstract Coxeter polyhedron is a \textit{potential hyperbolic Coxeter polyhedron}.\\
We say that a potential matrix $M$ \textit{specifies a potential compact hyperbolic Coxeter polyhedron} if each of its vertex groups (that is, the groups with Coxeter matrices of the form $M^s$, for $s$ a set of indices of facets intersecting at a vertex) is finite. We will refer to such a matrix as a \textit{potential matrix of compact type}.
\end{definition}

For each rank $r\geq 2$, there are infinitely many finite Coxeter groups, because of the infinite 1-parameter family of all dihedral groups, whose graphs consist of two nodes joined by an edge of weight $k\geq 2$. On the other hand, the following observation is immediate.
\begin{proposition}
There are finitely many finite Coxeter groups of rank $r$ with Coxeter matrix entries at most seven.
\end{proposition}
It thus suffices to enumerate matrices for potential Coxeter polyhedra with entries at most seven (and the remaining Coxeter matrices result from substituting integers greater than seven for pairs of sevens in these).\\

Since the Coxeter matrices of the irreducible finite and affine Coxeter groups are all known, determining the matrices for potential hyperbolic Coxeter polyhedra of a particular combinatorial type is a purely combinatorial problem. Not all such matrices will be Coxeter matrices of actually realizable hyperbolic Coxeter polyhedra of the desired combinatorial type, however. It is indeed necessary to also determine the length of the perpendicular segment between pairs of non-adjacent facets, and to check that the resulting Gram matrix has the right signature, namely $(n,1)$ (for more details on realizability conditions for hyperbolic Coxeter polyhedra, see \cite{Vinberg2}, for example).\\
The problem of finding such polyhedra is thus solved in two phases: first, we find potential matrices for hyperbolic Coxeter $n$-cubes, in particular using a backtracking search algorithm; and second, we solve the relevant algebraic conditions for the distances between non-adjacent facets (this will be the object of Section~\ref{sec grob}).

\subsection{Potential hyperbolic Coxeter $n$-cube matrices}\label{subsec matrices}
We consider now Coxeter $n$-cubes and study their potential matrices. Recall that we adopted the convention that the facets $2i-1$ and $2i$ of an $n$-cube do not intersect, $i=1,...,n$ (see Section \ref{subsec gencube}).

\begin{definition}
We call \textit{potential $n$-cube matrix} a potential matrix $M$ for an abstract polyhedron combinatorially equivalent to an $n$-cube.
\end{definition}

Let $M$ be a potential $n$-cube matrix. Then, in order for the associated abstract $n$-cube to be realizable as a hyperbolic polyhedron, the following necessary conditions must hold (see \cite[I, Chapter 6.2]{Vinberg2}, for example).
\begin{enumerate}
\item[(1)] For every sequence $s$ of indices of facets meeting at a vertex, the submatrix $M^s$ is the Coxeter matrix of a finite or affine Coxeter group.
\item[(2)] Any submatrix of $M$ of the form $M^{(2k-1,2k,2l-1,2l)}$ cannot be the Coxeter matrix of a Euclidean Coxeter 2-cube (see Remark~2.2 below).
\end{enumerate}

\begin{remark}
\begin{enumerate}
\item  Condition $(1)$ ensures that, if it is realizable, the associated hyperbolic polyhedron is really of finite volume (as it would be the convex hull of $2^n$ vertices in $\overline{\HH^n}$).
\item Condition $(2)$ is a special case of the so-called \textit{signature obstruction}: in the graph of a hyperbolic Coxeter polyhedron, any two subgraphs which are neither elliptic nor parabolic have to be connected by an edge (recall that the $\infty$ in positions $(2i-1,2i)$ of $M$ correspond to the dotted edges in the Coxeter graph).
\end{enumerate}
\end{remark}

Let $V:=\{1,2\}\times ...\times\{2n-1,2n\}$. Notice that by construction, the $2^n$ elements of $V$ correspond exactly to the respective labels of facets meeting at the vertices of an $n$-cube. To preserve the convention on labelling of opposite facets, consider the set $\Sigma$ of permutations $\sigma$ of $\{1,...,2n\}$ such that $\sigma$ induces a symmetry of the $n$-cube and such that for each $i=1,...,n$, there is a $j\in\{1,...,n\}$ such that $\{\sigma(2i-1),\sigma(2i)\}=\{2j-1,2j\}$.\\

Summing up conditions described above and in Section \ref{subsec gencube}, we get that a potential $n$-cube matrix is a symmetric $(2n)\times (2n)$ matrix $M$, with $1'$s on the diagonal and $\infty'$s in the odd-even positions of the superdiagonal and corresponding even-odd positions of the subdiagonal, having no Euclidean 2-cubes, and such that, for each $v\in V$, $M^v$ is the Coxeter matrix of a finite or affine Coxeter simplex group. We consider these matrices up to equivalence by relabelling, and so we only list those $M$ for which $M\leq M^\sigma$ for all $\sigma\in\Sigma$.\\

We order matrices lexicographically, taking $\infty$ before any number, by their lower triangular part, reading row-wise. Then, we enumerate these matrices with the the following \textit{backtracking search algorithm} (recall that it is sufficient to take $m_{ij}\in\{2,...,7\}\cup\{\infty\}$):
\begin{itemize}
\item An initial partial potential Coxeter $n$-cube matrix $M$ is taken, with 1's on the diagonal, $\infty$'s in the odd-even positions of the superdiagonal and corresponding even-odd positions of the subdiagonal, and with all remaining entries set to $\bigstar$.
\item At each stage, a vertex $n$-tuple $v\in V$ is taken.
\item The partial matrix $M^v$ is compared with the list of Coxeter matrices of finite and affine Coxeter simplex groups.
\item For each possible vertex Coxeter matrix refining $M^v$, trying each such matrix in turn, $\bigstar$ entries in $M$ are replaced by values to make $M^v$ the vertex Coxeter matrix.
\item If the partial matrix $M^v$ cannot be completed to the Coxeter matrix of a finite or affine Coxeter group, the algorithm backtracks (to the next alternative at an earlier stage).
\item If a Euclidean 2-cube is formed, the algorithm backtracks.
\item If a permutation of the partial matrix would be earlier in the ordering of matrices, then any completion would have an equivalent form earlier in the ordering and would have been enumerated by the algorithm in another sequence of choice, so the algorithm backtracks.
\item If every vertex has been checked (that is, if all submatrices $M^v$, $v\in V$, have been considered), then the partial matrix has been completed to a matrix satisfying all of the vertex constraints and which is earliest in the set of all equivalent forms.
\end{itemize}
At the end of the process, for each matrix with at least one pair of coefficients satisfying $m_{ij}=m_{ji}=7$ (corresponding to the presence of a dihedral subgroup in the associated Coxeter group, see Section~\ref{subsec potential}), remove from the list any matrix obtained by replacing $m_{ij}$ and $m_{ji}$ by $m\in\{2,...,6\}$. This removes redundant occurrences of Coxeter groups with dihedral subgroups.\\

This approach was programmed, and it successfully enumerated all potential Coxeter $n$-cube matrices for $n=3,4$ and $5$. However, a refined algorithm was needed to continue the search. The refined algorithm relies on the following observation.

\begin{lemma}\label{lemma subcube}
If $M$ is a potential Coxeter $n$-cube matrix, then $M^{(1,2,...,2n-2)}$ is a potential Coxeter $(n-1)$-cube matrix of compact type.
\begin{proof}
Let $M':=M^{(1,2,...,2n-2)}$. For each $(n-1)$-tuple $\sigma$ of $\{1,2\}\times...\times\{2n-3,2n-2\}$, the matrix $(M')^\sigma$ corresponds to a principal submatrix of the Coxeter matrix of a finite or affine Coxeter group. Such matrices are always Coxeter matrices of finite Coxeter groups. $_\blacksquare$
\end{proof}
\end{lemma}

For computational convenience, the ordering of $(2n)\times (2n)$ matrices entries is done first by their initial $(2n-2)\times(2n-2)$ submatrices, so that only an earliest representative of each $(n-1)$-dimensional matrix is needed to seed the search for the $n$-dimensional case. Then, the refined backtracking search is done one dimension at a time as follows.
\begin{itemize}
\item For each Coxeter $(n-1)$-cube matrix of compact type, start with a partial $(2n)\times(2n)$ matrix with the upper left $(2n-2)\times(2n-2)$ filled in.
\item Enumerate the possible $(2n-1)$-st rows, by imposing the vertex conditions involving facet $2n-1$ in a backtracking search as in the previous algorithm. 
\item Then, the $2n$-th row must satisfy the same conditions and there are no vertex conditions involving both facets $2n-1$ and $2n$ (since they are opposite in the cube). Hence, fill-in the rows $2n-1$ and $2n$ with each possible combination from the list of possible $(2n-1)$-st rows.
\item Remove from the list any matrices having earlier permutation equivalent form.
\end{itemize}

\subsection{Algorithm outcome}\label{subsec matrices outcome}
The above described procedure has been implemented in \textsf{Mathematica}. We describe the outcome, distinguishing between lower dimensions ($n=2,3$) and higher dimensions ($n\geq 4$).

\subsubsection{Dimensions 2 and 3}\label{ssubsec 23}
In dimension $n=2$, thank to a result of Poincar\'{e} (see \cite[I, Chapter 6.2]{Vinberg2}, for example), we know that the only necessary and sufficient condition for an abstract 2-cube to be realizable in $\mathbb{H}^2$ is that its four angles add up to less than $2\pi$. Hence, all potential hyperbolic Coxeter 2-cube matrices turn out to yield \textit{realizable} ones.\\
Recall that we enumerate matrices with entries at most seven (see Section \ref{subsec potential}). Replacing pairs of coefficients of the form $m_{ij}=m_{ji}=7$ with larger integers gives infinite families of potential Coxeter $n$-cubes.\\
The process begins with potential Coxeter 2-cube matrices of compact type, with entries at most seven. Up to relabeling, there are 230 such matrices, whose repartition is given in Table~\ref{table 2c} (see Theorem~\ref{thm classification}).\\
Notice that in this dimension, the only matrices corresponding to noncompact hyperbolic 2-cubes are exactly those containing at least one pair of coefficients of the form $m_{ij}=m_{ji}=\infty$.\\

In dimension 3, Andreev (see \cite[I, Chapter 6.2]{Vinberg2}, for example) gave a set of necessary and sufficient conditions, in terms of angles and combinatorics, for an abstract 3-polyhedron to be realizable in the hyperbolic space. In our context, the combinatorial type is fixed, and Andreev's Theorem implies that every potential hyperbolic Coxeter 3-cube matrix is in fact realizable. Andreev's result, however, does not provide the respective lengths of the perpendicular segments between non-intersecting facets, and, more importantly, it does not give any clue about the size of the family being considered.\\
As a result, we get the repartition given in Table~\ref{table 3c} (see Theorem~\ref{thm classification}). In particular, the number in each cell of Table~\ref{table 3c} indicates the number of hyperbolic Coxeter $3$-cubes with no free parameter (that is, with angles at most $\frac{\pi}{6}$), respectively the number of infinite families of hyperbolic Coxeter $n$-cubes indexed by $k$ free parameters, $k=1,2,3$, corresponding to dihedral subgroups in the associated Coxeter group (see Section~\ref{subsec potential}).\\
In particular, there is no family of hyperbolic Coxeter $3$-cubes with 4 or more free parameters.

\begin{remark}
\begin{itemize}
\item We provide the explicit list of all hyperbolic Coxeter 3-cubes as an Appendix (see p.22).
\item The classification described in Table~\ref{table 3c} contains several already available partial classifications (see \cite{CHL,ImHof,Jacquemet3}, for example).
\item The unique 3-parameter family of hyperbolic Coxeter 3-cubes is the well-known family of \textit{Lambert cubes} \cite{ImHof}.
\item The large amount of single occurrences and the presence of infinite families, both in the compact and noncompact settings, contrast with the simplex case, where the classification is reduced to 9 compact and 22 noncompact hyperbolic Coxeter 3-simplices \cite[pp 205-208]{Vinberg2}.
\end{itemize}
\end{remark}

\subsubsection{Higher dimensions}\label{ssubsec higher}
The repartition of potential hyperbolic Coxeter $n$-cubes matrices, $n=4,5,6$, is given in Table~\ref{table pm}.

\begin{table}[h!]
\centering \renewcommand{\arraystretch}{1.25}
\begin{tabular}{| c | c|| c | c | c |}
\hline Type & $k$ & $n=4$ & $n=5$ & $n=6$\\
\hline\hline
\multirow{3}*{Compact} & 0 & 7419 & 589 & 0\\
\cline{2-5}
& 1 & 188 & 1 &0 \\
\cline{2-5}
& $\geq 2$ & 0 & 0 & 0\\
\hline
\multirow{4}*{Noncompact} & 0 & 19068 & 3111 & 2\\
\cline{2-5}
& 1 & 671 & 18 & 0 \\
\cline{2-5}
& 2 & 21 & 0 & 0\\
\cline{2-5}
& $\geq 3$ & 0 & 0 & 0\\
\hline
\end{tabular}
 \caption{Number of potential hyperbolic Coxeter $n$-cube matrices with $k$ free parameters, $n=4,5,6$}\label{table pm}
\end{table}

This leads immediately to the following observation.

\begin{corollary}\label{cor high}
There are no hyperbolic Coxeter $n$-cubes for any $n\geq 7$.
\end{corollary}
\begin{proof}
This follows directly from Lemma~\ref{lemma subcube} and the fact that there are no potential hyperbolic Coxeter $6$-cube matrices of compact type (see Table~\ref{table pm}). $_\blacksquare$
\end{proof}

\section{From potential matrices to realizable matrices}\label{sec grob}
It remains now to check each potential Coxeter $n$-cube matrix to see if it can lead to the Gram matrix $G$ of an actual hyperbolic Coxeter $n$-cube. 
Recall that for $n=2,3$, any potential hyperbolic Coxeter $n$-cube matrix is automatically realizable (see Section \ref{ssubsec 23}), so that the classification problem there reduces to the computation of the respective length between pairs of non-intersecting facets. Hence, it remains to consider the dimensions $n=4,5,6$.\\
In the sequel, $M$ will denote one of the potential hyperbolic Coxeter $n$-cube matrices obtained in Section \ref{ssubsec higher}. We follow the same labelling conventions as in Section \ref{subsec matrices}.

\subsection{Algebraic realizability conditions}\label{subsec conditions}
The $n$ perpendicular line segments between pairs of non-adjacent facets of an $n$-cube $\mathcal{P}$ determine the entries of the Gram matrix corresponding to the $\infty$ entries at positions $(2i-1,2i)$ and $(2i,2i-1)$ in the potential Coxeter matrix $M$. The corresponding Gram matrix entries will be of the form $-\cosh(d_i)$, for $d_i$ the distance between the facets $2i-1$ and $2i$ of $\mathcal{P}$.\\
The condition on these entries is that the resulting Gram matrix must have signature $(n,1)$. This implies that the following set of polynomial conditions must be satisfied:
\begin{equation}\label{eq detG}
\mathrm{det}(G^S)=0 \text{ for any } (n+2)\text{-tuple } S \text{ of facet indices}.
\end{equation}
Thus, a first natural idea is to take the Gram matrix $G$ with entries $-\cos(\pi/m_{ij})$ for a potential matrix $M=(m_{ij})_{1\leq i,j\leq 2n}$, replacing the entries $(2i-1,2i)$ and $(2i,2i-1)$ by a variable $x_i$ for each $i=1,...,n$, to solve Conditions $(\ref{eq detG})$, with the constraints
\begin{equation}\label{eq xi}
x_i<-1,\,\,\, i=1,...,2n,
\end{equation}
and finally to check that the resulting Gram matrix has the right signature~$(n,1)$.\\

Unfortunately, this is not quite so easy (even for a computer), in view of the mass production scale of the numbers of potential Coxeter $n$-cube matrices we have to consider. One difficulty has to do with the infinite families of potential Coxeter matrices we have identified, so that each pair of entries of the form $m_{ij}=m_{ji}=7$ can be replaced by values greater than 7. To alleviate this issue, we refine our approach by adding further variables $z_{\{i,j\}}$ corresponding to each such pair of sevens in the potential matrix, and replace the entries $-\cos(\pi/m_{ij})$ and $-\cos(\pi/m_{ji})$ of $G$ by $-z_{\{i,j\}}/2$, so that $z_{\{i,j\}}=2\cos(\pi/m_{ij})$. Then, the associated constraint is 
\begin{equation}\label{eq zij}
2\cos(\pi/7)\leq z_{\{i,j\}}<2,\text{ for all }\{i,j\}\text{ such that }m_{ij}\geq 7.
\end{equation}

Hence, by Table~\ref{table pm} and Conditions (\ref{eq detG}), (\ref{eq xi}) and (\ref{eq zij}), we have the following respective polynomial conditions.
\begin{itemize}
\item In dimension 4, there are 24 polynomial conditions to consider, in four $x_i$ variables and at most two $z_{\{i,j\}}$ variables, each with their own inequality constraint.
\item In dimension 5, there are 80 polynomial conditions in five $x_i$ variables and at most one $z_{\{i,j\}}$ variable.
\item In dimension 6, there are 240 polynomial conditions with six $x_i$ variables only.
\end{itemize}
The main problem is to efficiently eliminate from consideration all the non-realizable potential Coxeter $n$-cubes matrices. As a technical contrivance, it is convenient to introduce variables $t_4$, $t_5$, $t_6$ and substitute any Gram matrix entries of the form $-\cos(\pi/m)$ with $-t_m/2$, for $m=4,5$ or $6$, with then additional equations
\begin{equation}\label{eq tieq}
t_4^2-2=0,\,\,\,t_5^2-t_5-1=0 \text{ and }t_6^2-3=0,
\end{equation}
and additional inequalities
\begin{equation}\label{eq tiineq}
t_m>0,\,\,\,m=4,5,6.
\end{equation}
Then, the conditions 
\begin{equation}\label{eq gram2}
\mathrm{det}(2G^S)=0 \text{ for any } (n+2)\text{-tuple } S \text{ of facet indices}.
\end{equation}
are equivalent to Conditions~(\ref{eq detG}) and are integer polynomial equations in the $t_m$'s, $x_i$'s and $z_{\{i,j\}}$'s.

\begin{remark}
Since there are so many conditions and so few variables to determine, it is to expect that not many of the potential Coxeter $n$-cube matrices will give rise to hyperbolic Coxeter $n$-cubes, as we will see.
\end{remark}

A powerful systematic approach is to find \textit{Gröbner bases} for the sets of polynomials arising from $(\ref{eq tieq})$ and $(\ref{eq gram2})$. We briefly outline this concept and give some useful references.\\
Let $K\subset\mathbb{C}$ be a field, and let $\mathcal{F}$ be a finite set of polynomials in $K[z_1,...,z_l]$, $l\geq 1$ (here, $K=\mathbb{R}$ and $\mathcal{F}=\{\mathrm{det}(2G^S)\,|\,S \text{ is a } (n+2)-\text{tuple of facet indices}\} \cup\{t_4^2-2,\,t_5^2-t_5-1,\,t_6^2-3\}$). Recall that a Gröbner basis $\mathcal{G}$ for $\mathcal{F}$ is a finite set of polynomials enjoying particularly nice properties (such as having 'easy to determine' zeroes), and such that the variety of $\mathcal{F}$ coincides with the variety of $\mathcal{G}$. In particular, one has $\mathcal{G}=\{1\}$ if and only if the polynomials in $\mathcal{F}$ do not share common zeroes.\\
Formal definitions as well as key results and algorithms to extract Gröbner bases from set of polynomials can be found in \cite{CLO,VzG}, for example. This procedure is well understood, and can be performed by almost any mathematical software. For this work, we have used the \textsf{Mathematica} packages.

\subsection{Gröbner bases extraction}\label{subsec grobner}

\subsubsection{Dimension 6}
There are only 2 potential matrices in dimension 6, giving rise to the following candidates for Gram matrices of hyperbolic Coxeter $6$-cubes:
{\tiny \[G_1:=\left(
\begin{array}{cccccccccccc}
 1 & x_1 & 0 & 0 & 0 & 0 & 0 & -1/2 & 0 & 0 & 0 & -1/2 \\
 x_1 & 1 & 0 & -1/2 & 0 & -1/2 & 0 & 0 & 0 & -1/2 & 0 & 0 \\
 0 & 0 & 1 & x_2 & 0 & -1/2 & 0 & 0 & 0 & 0 & -1/2 & 0 \\
 0 & -1/2 & x_2 & 1 & 0 & 0 & 0 & -1/2 & -1/2 & 0 & 0 & 0 \\
 0 & 0 & 0 & 0 & 1 & x_3 & -1/2 & 0 & 0 & -1/2 & 0 & 0 \\
 0 & -1/2 & -1/2 & 0 & x_3 & 1 & 0 & 0 & 0 & 0 & 0 & -1/2 \\
 0 & 0 & 0 & 0 & -1/2 & 0 & 1 & x_4 & -1/2 & 0 & 0 & -1/2 \\
 -1/2 & 0 & 0 & -1/2 & 0 & 0 & x_4 & 1 & 0 & 0 & 0 & 0 \\
 0 & 0 & 0 & -1/2 & 0 & 0 & -1/2 & 0 & 1 & x_5 & -1/2 & 0 \\
 0 & -1/2 & 0 & 0 & -1/2 & 0 & 0 & 0 & x_5 & 1 & 0 & 0 \\
 0 & 0 & -1/2 & 0 & 0 & 0 & 0 & 0 & -1/2 & 0 & 1 & x_6 \\
 -1/2 & 0 & 0 & 0 & 0 & -1/2 & -1/2 & 0 & 0 & 0 & x_6 & 1 \\
\end{array}
\right),\]}

and

{\tiny \[ G_2:=\left(
\begin{array}{cccccccccccc}
 1 & x_1 & 0 & 0 & 0 & 0 & 0 & -1/2 & 0 & -1/2 & 0 & 0 \\
 x_2 & 1 & 0 & -1/2 & 0 & -1/2 & 0 & 0 & 0 & 0 & 0 & -1/2 \\
 0 & 0 & 1 & x_2 & 0 & -1/2 & 0 & 0 & 0 & -1/2 & -1/2 & 0 \\
 0 & -1/2 & x_2 & 1 & 0 & 0 & 0 & -1/2 & 0 & 0 & 0 & 0 \\
 0 & 0 & 0 & 0 & 1 & x_3 & 0 & -1/2 & -1/2 & 0 & -1/2 & 0 \\
 0 & -1/2 & -1/2 & 0 & x_3 & 1 & 0 & 0 & 0 & 0 & 0 & 0 \\
 0 & 0 & 0 & 0 & 0 & 0 & 1 & x_4 & 0 & -1/2 & 0 & -1/2 \\
 -1/2 & 0 & 0 & -1/2 & -1/2 & 0 & x_4 & 1 & 0 & 0 & 0 & 0 \\
 0 & 0 & 0 & 0 & -1/2 & 0 & 0 & 0 & 1 & x_5 & 0 & -1/2 \\
 -1/2 & 0 & -1/2 & 0 & 0 & 0 & -1/2 & 0 & x_5 & 1 & 0 & 0 \\
 0 & 0 & -1/2 & 0 & -1/2 & 0 & 0 & 0 & 0 & 0 & 1 & x_6 \\
 0 & -1/2 & 0 & 0 & 0 & 0 & -1/2 & 0 & -1/2 & 0 & x_6 & 1 \\
\end{array}
\right).\]}

Notice that there is no $t_m$ variable here. The Gröbner basis for the polynomials arising from Conditions (\ref{eq gram2}) for $G_1$ can be computed to be $\{1\}$, so $G_1$ cannot be completed to a Gram matrix of an actual hyperbolic Coxeter $6$-cube.\\
The Gröbner basis arising form $G_2$ includes the polynomials $2x_i-1$, $i=1,...,6$, so that the only solution of the polynomial conditions is $x_i=1/2$, $i=1,...,6$, and this is inconsistent with Condition (\ref{eq xi}). \\
This proves the following result, which is part 4 of Theorem~\ref{thm classification}.
\begin{proposition}\label{prop 6cubes}
There is no hyperbolic Coxeter 6-cube.
\end{proposition}

\begin{remark}
As it is a negative result, Proposition~\ref{prop 6cubes} can already be obtained by exhibiting an appropriate impossible setting. For instance, for $G_2$, it is sufficient to see that Condition~(\ref{eq gram2}) for $S=\{1,...,12\}\setminus\{2,5,10,12\}$ leads to the condition
\[1 - 2 x_4^2 - 3 x_2^2 + 4 x_2^2 x_4^2=0,\]
which cannot be satisfied for $x_2,x_4<-1$.\\
However, this exclusion strategy (which has been successfully applied in the particular case treated in \cite{Jacquemet3}), in contrast to the Gröbner basis extraction algorithm, cannot be made systematic, as there is no satisfactory method to decide how the sets $S$ of indices have to be chosen, and how many of such sets are sufficient to yield an impossible setting.
\end{remark}

\subsubsection{Dimension 5}\label{subsec dim5}
In dimension 5, it took \textsf{Mathematica} about 10 minutes to compute the Gröbner bases for the 3719 potential matrices, of which 3707 reduce to $\{1\}$.\\
Of the remaining 12 cases, one basis includes the polynomials $2x_i-1$, $i=1,...,5$, and another one includes the polynomial $x_5-1$, each eliminated by the conditions $x_i<-1$, $i=1,...,5$.\\
Another 6 of the remaining potential Coxeter 5-cubes matrices have a pair of seven entries introducing a variable $z_{\{i,j\}}=:z$. Each of these cases has a solution satisfying $x_i<-1$, $i=1,...,5$, but with $z$ outside the range $[2\cos(\pi/7),2)$ given by Condition~(\ref{eq zij}). In fact, one of these cases had $z=0$, corresponding to $m=2$, and substituting two for the sevens in the potential Coxeter matrix gave (a permuted form of) one of the remaining four cases, while another of these cases had $z=1$, corresponding to $m=3$, and substituting threes for the sevens in the potential Coxeter matrix gave another (permuted form of) one of the remaining four potential matrices, so that these two cases can be considered as redundant.\\
The remaining four potential matrices all lead to solutions satisfying all of the conditions and having signature $(5,1)$, hence proving Part 3 of Theorem~\ref{thm classification}.\\

The respective lengths between pairs of non-intersecting facets are given as follows. Let $\cosh\,d$ denote the (equal) weights of the two dotted edges lying on the diagonal of the cube formed by the graph $\Gamma_1$ of the compact hyperbolic Coxeter $5$-cube (see Figure~\ref{fig 5cubec}), and let $\cosh\,e$ denote the (equal) weights of the remaining three dotted edges. Then, one has $\cosh\,d=\frac{\sqrt{2(5+\sqrt{13})}}{4}$, and $\cosh\,e=\frac{1+\sqrt{13}}{4}$.\\

As for the \text{graphs of noncompact type} (see Figure~\ref{fig 5cubesnc}), we have the following situation. Let $\cosh\,d^m_{i}$ be the weight associated to the dotted edge connecting the vertices $2i-1$ and $2i$ in $\Gamma^m_2$, $i=1,...,5$, $m=2,3$. Then, the respective weights are given by
\[\cosh\,d^m_i=\frac{\sqrt{2}+\sqrt{10}}{4},\,i=1,...,4,\,m=2,3,\]
and
\[\cosh\,d^2_5=\frac{\sqrt{5}}{2},\text{ and }\cosh\,d_5^3=\frac{\sqrt{2+\sqrt{5}}}{2}.\]

Finally, the weights of the dotted edges of $\Gamma_3$ are all the same, and they are given by $\cosh\,d=\frac{\sqrt{5}}{2}$.

\begin{remark}
By Vinberg's arithmeticity criterion \cite[Section 5]{Vinberg67}, one can deduce that the groups with Coxeter graph $\Gamma_1$ and $\Gamma_3$ are arithmetic, and that the groups with Coxeter graph $\Gamma^m_2$, $m=2,3$, are non-arithmetic.
\end{remark}

\begin{remark}\label{rem com}
Although the commensurability problem is difficult in general, available tools (see \cite{GJK, JKRT2}, for example) can be used to show that the \textit{arithmetic non-cocompact} Coxeter group with Coxeter graph $\Gamma_3$ is commensurable to the Coxeter simplex group with Coxeter symbol $[3,3^{[5]}]$. 
\end{remark}

\subsubsection{Dimension 4}\label{subsec dim4}
The efficiency of the process for finding a Gröbner basis, even when the result is just $\{1\}$, depends significantly on the ordering of the variables, leading to possible technical difficulties: intermediate steps in the Gröbner basis extraction procedure may involve polynomials of high degree with huge coefficients. To make efficient progress in the dimension 4 case, it was necessary to run the Gröbner basis routine for different subsets of equations and variables with time limits on execution, aborting slow calculations in hopes of finding a method of completing the calculation in a reasonable time.\\

We provide some details on the computation procedures and timing limitations, as a hint on the difficulties which can be expected during the Gröbner bases extraction procedure. Notice that the exact number of cases resulting on a program time out may vary somewhat between runs (depending on the computer and what other tasks are competing for computing time), so that the numbers provided here should be considered as an order of magnitude on the complexity of the procedure.
\begin{enumerate}
\item The equations for the 27367 potential matrices took about 15 minutes to be generated by \textsf{Mathematica}. An initial pass over all cases using all equations and allowing a maximum of 5 seconds per case finds a Gröbner basis in 26967 cases in about 75 minutes, and all but 1214 of these reduce to $\{1\}$ 
\item A second pass over the remaining 400, first finding a Gröbner basis for the equations not involving $t_5$ and $x_1$ before reducing with the remaining equations, and allowing 10 seconds per case, took about 7 minutes and found another 368 Gröbner bases, all reducing to $\{1\}$.
\item A third pass over the remaining 32 cases, first finding a Gröbner basis for the equations not involving $t_4$ and $x_2$ before reducing with the remaining equations, and allowing 20 seconds per case, took less than a minute and found the remaining 32 Gröbner bases, all reducing to $\{1\}$.
\end{enumerate}

Hence, it remains to analyse the 1214 potential matrices leading to Gröbner bases different from $\{1\}$ obtained in Step 1 above. These can be considered in 14 separate cases depending on the set of values not belonging to $\{1,2,3,\infty\}$, and the number of pairs of sevens (each possibly substituted with a greater label), in the Coxeter matrix. In other words, the cases are described as depending on occurrences of $t_4-\sqrt{2}$, $t_5-(1+\sqrt{5})/2$, or $t_6-\sqrt{3}$ in the Gröbner basis, and on the number of additional $z_{\{i,j\}}$ variables, respectively.\\
These 14 cases and the respective numbers of corresponding potential matrices are as follows (recall that from Table~\ref{table pm}, it follows that no potential hyperbolic Coxeter 4-cube matrix has more than 2 pairs of sevens).
\begin{itemize}
\item With two pairs of sevens, there are:
\begin{enumerate}
\item[(1)] 6 cases with no 4, 5, or 6 label, and
\item[(2)] 15 cases including a 4 label and no 5 or 6 label.
\end{enumerate}
\item With a single pair of sevens, there are:
\begin{enumerate}[resume]
\item[(3)] 31 cases with no 4, 5, or 6 label,
\item[(4)] 223 cases with a 4 label but no 5 or 6 label,
\item[(5)] 18 cases with a 5 label but no 4 or 6 label,
\item[(6)] 9 cases with a 6 label but no 4 or 5 label,
\item[(7)] 47 cases with a 4 and 5 label but no 6 label, and
\item[(8)] 27 cases with a 4 and 6 label but no 5 label. 
\end{enumerate}
\item With no sevens, there are:
\begin{enumerate}[resume]
\item[(9)] 56 cases with no 4, 5 or 6 label, 
\item[(10)] 538 cases with a 4 label but no 5 or 6 label,
\item[(11)] 59 cases with a 5 label but no 4 or 6 label,
\item[(12)] 3 cases with a 6 label but no 4 or 5 label,
\item[(13)] 168 cases with a 4 and 5 label but no 6 label, and
\item[(14)] 14 cases with a 4 and 6 label but no 5 label. 
\end{enumerate}
\end{itemize}

Then, we proceed by setting each polynomial in the Gröbner bases to zero, with the inequality constraints $x_i<-1$ and $2\cos(\pi/7)<z_{\{i,j\}}<2$ for any added $z_{\{i,j\}}$ variable, in order to solve the corresponding system of polynomial equations. As for the Gröbner bases extraction routine, these reductions were performed with time constraints and recursively performed for different orders of the variables and starting with different subsets of equations, since again the intermediate stages of reduction could result in large expressions even if the final result was no solution. This also allowed us to keep track of the procedure progress.\\
As a result, cases (1) to (8) could be exhausted (in less than 25 minutes by \textsf{Mathematica}), yielding no solutions in either case, and all cases (9)-(14) were treated efficiently in less than a minute combined (due to the absence of $z_{\{i,j\}}$ variable), returning exactly 15 \textit{realizable} potential matrices.\\
This proves Part 2 of Theorem~\ref{thm classification}, and therefore achieves its proof.\\

The respective weights of the dotted edges in the graphs $\Sigma_1^{k,l,m}$ and $\Sigma_2^m$ (see Figure~\ref{fig 4cubes}) are given in Tables~\ref{table dottedc41} and \ref{table dottedc42}, respectively.\\
\begin{table}[h!]
\centering \renewcommand{\arraystretch}{1.9}
\begin{tabular}{| c | c || c | c | c |}
\hline \multicolumn{2}{|c||}{}  & $\Sigma^{k,l,3}_1$ & $\Sigma^{k,l,4}_1$ & $\Sigma^{k,l,5}_1$\\
\hline\hline
\multirow{3}*{$\cosh\,d$} & $(k,l)=(2,2)$ & $\sqrt{\frac{5+\sqrt{13}}{6}}$ & $\frac{\sqrt{2}+\sqrt{10}}{4}$ & $\sqrt{\frac{3}{2}-\frac{\sqrt{5}}{10}}$\\
\cline{2-5}
& $(k,l)=(2,3)$ & $\frac{2+\sqrt{13}}{3}$ & $\frac{1+\sqrt{5}}{2}$ & $2-\frac{\sqrt{5}}{5}$ \\
\cline{2-5}
& $(k,l)=(3,2)$ & $\frac{5+\sqrt{13}}{4}$ & $\frac{6+\sqrt{2}+2\sqrt{5}+\sqrt{10}}{8}$ & $\frac{5+\sqrt{5}}{4}$\\
\hline

\multicolumn{2}{|c||}{$\cosh\,e$} &  $\frac{1+\sqrt{13}}{4}$ & $\frac{\sqrt{2}+\sqrt{10}}{4}$ & $\frac{\sqrt{5}}{2}$\\
\hline
\end{tabular}
\caption{Weights of the dotted edges in the graphs $\Sigma^{k,l,m}_1$. Here, $\cosh\,d$ refers to the weight of the dotted edge which is the diagonal of the graph of $\Sigma^{k,l,m}_1$, and $\cosh\,e$ refers to the weight of the three dotted edges which are (boundary) edges of the graph of $\Sigma^{k,l,m}_1$ (in particular, these three dotted edges always have the same weight, which does not depend on the parameters $k$ and $l$).}\label{table dottedc41}
\end{table}

\begin{table}[h!]
\centering \renewcommand{\arraystretch}{1.9}
\begin{tabular}{| c || c | c |}
\hline  & $\Sigma^2_2$ & $\Sigma^3_2$ \\
\hline\hline
$\cosh\,d_1$ & $\frac{1+\sqrt{5}}{2}$ & $\frac{1+\sqrt{5}}{2}$\\
\hline
$\cosh\,d_2$ & $\frac{1+\sqrt{5}}{2}$ & $\frac{\sqrt{2}+\sqrt{10}}{4}$ \\
\hline
$\cosh\,d_3$ & $\frac{\sqrt{2}+\sqrt{10}}{4}$  & $\frac{\sqrt{2}+\sqrt{10}}{4}$ \\
\hline
\end{tabular}
 \caption{Weights of the dotted edges in the graphs $\Sigma_2^{m}$, $m=2,3$. Here, $\cosh\,d_1$ refers to the weight of the horizontal dotted edge in the graph of $\Sigma^m_2$, $\cosh\,d_2$ to the weight of the vertical dotted edge, and $\cosh\,d_3$ to the weight of both diagonal dotted edges in the graph of $\Sigma_2^m$.}\label{table dottedc42}
\end{table}

As for the graph $\Sigma_3$, all its dotted edges have the same weight, which is given by $\cosh\,d=\frac{1+\sqrt{2}}{2}$.

\begin{remark}
The volumes of the hyperbolic Coxeter 4-cubes can be computed, for example by using the computer program \textsf{CoxIter} \cite{Raf}. They are given in Table~\ref{table volc4}. Notice that for these polyhedra, the relations of the form $\text{vol}(P_1)=2\,\text{vol}(P_2)$ result from the fact that $P_1$
is a suitable doubling of $P_2$.
\end{remark}
\begin{table}[H]
 \begin{center}\renewcommand{\arraystretch}{1.3}
 \begin{tabular}{| m{0.14\linewidth} || m{0.09\linewidth} | m{0.09\linewidth} | m{0.09\linewidth} | m{0.09\linewidth} | m{0.09\linewidth} |}
    \hline
    \centering $\text{graph }\Sigma(\mathcal{P})$
    & \centering $\Sigma^{2,2,3}_1$ 
    & \centering $\Sigma^{2,2,4}_1$
	& \centering $\Sigma^{2,2,5}_1$
	& \centering $\Sigma^2_2$
	& \centering $\Sigma_3$
	\tabularnewline \hline \hline
	\centering $\mathrm{vol}\,\mathcal{P}$\rule[-0.35cm]{0cm}{0.9cm}
    & \centering $\frac{29}{720}\pi^2$ 
    & \centering $\frac{17}{288}\pi^2$ 
    & \centering  $\frac{133}{1800}\pi^2 $
    & \centering $\frac{17}{216}\pi^2 $
    & \centering $\frac{131}{1920}\pi^2 $
	\tabularnewline \hline
 \end{tabular}
  \caption{Volumes of some hyperbolic Coxeter $4$-cubes. The remaining volumes are given according to the relations $\text{vol}\,\mathcal{P}_1^{2,3,m}=\text{vol}\,\mathcal{P}^{3,2,m}_1=2\,\text{vol}\,\mathcal{P}^{2,2,m}_1$, $m=3,4,5$, and $\text{vol}\,\mathcal{P}^3_2=2\,\text{vol}\,\mathcal{P}^2_2$, where $\mathcal{P}_i^j$ is the 4-cube with graph $\Sigma(\mathcal{P}_i^j)=\Sigma_i^j$.}\label{table volc4}
\end{center}
\end{table}

Finally, for the graphs of noncompact type $\Sigma_4^m$, $m=2,3$, and $\Sigma_5$ (see Figure~\ref{fig 5cubesnc}), we have the following situation.
\begin{itemize}
\item Let $\cosh\,d_m$ denote the weight of the horizontal dotted edge in $\Gamma_4^m$, and let $\cosh\,e_m$ denote the (coinciding) weights of the vertical dotted edges in $\Gamma_4^m$. Then, we have $\cosh\,d_2=3/2=\cosh\,d_3$, $\cosh\,e_2=2\sqrt{2/7}$, and $\cosh\,e_3=9/7$.\\
Moreover, all the dotted edges in $\Sigma_5$ have the same weight, which is given by $\cosh\,d=3/2$.
\item For the volumes, one has $2\,\text{vol}\,\mathcal{P}_4^2=\text{vol}\,\mathcal{P}_4^3=\text{vol}\,\mathcal{P}_5=2\pi^2/9$.
\end{itemize}

\begin{remark}
As for the 5-dimensional case, Vinberg's arithmeticity criterion \cite[Section 5]{Vinberg67} allows us to conclude that the groups with Coxeter graph $\Sigma_1^{2,2,3}$, $\Sigma_1^{k,l,4}$, $\Sigma_1^{k,l,5}$, and $\Sigma^m_2$ are non-arithmetic and not quasi-arithmetic, that the groups with Coxeter graph $\Sigma^{2,3,3}_1$ and $\Sigma_4^m$ are quasi-arithmetic, and that the groups with Coxeter graph $\Sigma_1^{3,2,3}$, $\Sigma_3$ and $\Sigma_5$ are arithmetic.
\end{remark}

\begin{remark}
Similarly to the 5-dimensional case (see Remark~\ref{rem com}), one can show that the \textit{arithmetic non-cocompact} Coxeter group with Coxeter graph $\Sigma_5$ is commensurable to the Coxeter simplex group with Coxeter symbol $[3,4,3,4]$. 
\end{remark}

\begin{center}\subsubsection*{Acknowledgements}\end{center}
The first author has been supported by the Swiss National Science Foundation, fellowship Nr. P2FRP2-161727.\\
The authors would like to thank Rafael Guglielmetti for his comments and his assistance in computing invariants of the hyperbolic Coxeter 4- and 5-cubes, Ruth Kellerhals and John Ratcliffe for their useful comments on a preliminary version of this paper, and the anonymous referees for their helpful comments and suggestions.\\

\bibliographystyle{abbrv}
\bibliography{biblio}

\section*{Appendix - All hyperbolic Coxeter 3-cubes}
In the sequel, a hyperbolic Coxeter 3-cube is denoted by the 12-tuple of Coxeter matrix entries below the diagonal coefficients $m_{ii}=1$, $i=1,...,6$, and
subdiagonal coefficients $m_{21}=m_{43}=m_{65}=\infty$ listed in order
$$(m_{31},m_{32},m_{41},m_{42},\ \ m_{51},m_{52},m_{53},m_{54},\ \ m_{61},m_{62},m_{63},m_{64}).$$
An entry $m_{ij}=7$ can be replaced by any higher integer value in an infinite family. Only the Coxeter matrices having lexicographically first
12-tuple is given (any other matrix would be equivalent to one of the matrices listed here, by permuting rows and column in the right way).\\
The first four entries determine the 2-cube ``seed'' of the 3-cube (see Lemma~\ref{lemma subcube}, p.11). It turns out that of the 230 compact cases in dimension 2, only 26 occur as such a (lexicographically earliest) seed.\\
Given the seed, the next four values determine the angles to adjacent facets for facet 5, constrained to make four of the vertex groups
finite or affine. But then the exact same conditions apply for the last four values and facet 6. The 4-tuple for facet 5 will be lexicographically
the same or before that for facet 6, since otherwise the symmetry swapping facets 5 and 6 would give a lexicographically earlier 12-tuple.\\

We give tables for each $(m_{31},m_{32},m_{41},m_{42})$, each table with columns for each $(m_{51},m_{52},m_{53},m_{54})$ column labels
(writtem vertically), and then rows for each (equal or later) $(m_{61},m_{62},m_{63},m_{64})$ marking the table entry if there is
a 3-cube and with this as lexicographically specification by a ``c'' for compact and ``n'' for noncompact, and 
otherwise taking a period ``.'' entry.

\begin{table}[h!]
{\centering\setlength{\tabcolsep}{0.44pt}\footnotesize
\begin{tabular}{cccccccccccccccccccccccccccccccccccccccccccccccc}
 $m_{31}$ & $m_{32}$ & $m_{41}$ & $m_{42}$ & $m_{51}$ & 2 & 2 & 2 & 2 & 2 & 2 & 2 & 2 & 2 & 2 & 2 & 2 & 2 & 2 & 2 & 2 & 2 & 2 & 2 & 2 & 2 & 2 & 2 & 2 & 2 & 2 & 2 & 2 & 2 & 2 & 2 & 2 & 2 & 2 & 2 & 2 & 2 & 2 & 2 & 2 & 2 & 2 & 2 \\
 2 & 2 & 2 & 3 & $m_{52}$ & 2 & 2 & 2 & 2 & 2 & 2 & 2 & 2 & 2 & 2 & 2 & 2 & 2 & 2 & 2 & 2 & 2 & 2 & 2 & 2 & 2 & 2 & 2 & 2 & 2 & 2 & 2 & 2 & 2 & 2 & 3 & 3 & 3 & 3 & 3 & 3 & 3 & 3 & 3 & 4 & 4 & 5 & 6 \\
   &   &   &   & $m_{53}$ & 2 & 2 & 2 & 2 & 2 & 3 & 3 & 3 & 3 & 3 & 4 & 4 & 4 & 4 & 4 & 5 & 5 & 5 & 5 & 5 & 6 & 6 & 6 & 6 & 6 & 7 & 7 & 7 & 7 & 7 & 2 & 3 & 3 & 4 & 4 & 5 & 5 & 6 & 6 & 3 & 4 & 3 & 3 \\
 $m_{61}$ & $m_{62}$ & $m_{63}$ & $m_{64}$ & $m_{54}$ & 2 & 3 & 4 & 5 & 6 & 2 & 3 & 4 & 5 & 6 & 2 & 3 & 4 & 5 & 6 & 2 & 3 & 4 & 5 & 6 & 2 & 3 & 4 & 5 & 6 & 2 & 3 & 4 & 5 & 6 & 3 & 2 & 3 & 2 & 3 & 2 & 3 & 2 & 3 & 2 & 2 & 2 & 2 \\
 2 & 3 & 2 & 3 &   & n & . & . & . & . & n & . & . & . & . & n & . & . & . & . & n & . & . & . & . & n & . & . & . & . & n & . & . & . & . & n & . & . & . & . & . & . & . & . & . & . & . & . \\
 2 & 3 & 3 & 2 &   & c & c & . & . & . & . & c & . & . & . & . & c & . & . & . & . & c & . & . & . & . & c & . & . & . & . & c & . & . & . & n & c & . & . & . & . & . & . & . & . & . & . & . \\
 2 & 3 & 3 & 3 &   & n & n & . & . & . & n & n & . & . & . & n & n & . & . & . & n & n & . & . & . & n & n & . & . & . & n & n & . & . & . & n & n & n & . & . & . & . & . & . & . & . & . & . \\
 2 & 3 & 4 & 2 &   & c & c & c & . & . & . & c & c & . & . & . & c & c & . & . & . & c & c & . & . & . & c & c & . & . & . & c & c & . & . & n & c & n & c & . & . & . & . & . & . & . & . & . \\
\end{tabular}
\par}
\caption{Cubes $(2,2,2,3,\ \ 2,\ast,\ast,\ast,\ \ 2,\ast,\ast,\ast)$}
\end{table}

\begin{table}[h!]
{\centering\setlength{\tabcolsep}{0.44pt}\footnotesize
\begin{tabular}{cccccccccccccccccccccccccccccccccccccccccccccccc}
 $m_{31}$ & $m_{32}$ & $m_{41}$ & $m_{42}$ & $m_{51}$ & 2 & 2 & 2 & 2 & 2 & 2 & 2 & 2 & 2 & 2 & 2 & 2 & 2 & 2 & 2 & 2 & 2 & 2 & 2 & 2 & 2 & 2 & 2 & 2 & 2 & 2 & 2 & 2 & 2 & 2 & 2 & 2 & 2 & 2 & 2 & 2 & 2 & 2 & 2 & 2 & 2 & 2 & 2 \\
 2 & 2 & 2 & 3 & $m_{52}$ & 2 & 2 & 2 & 2 & 2 & 2 & 2 & 2 & 2 & 2 & 2 & 2 & 2 & 2 & 2 & 2 & 2 & 2 & 2 & 2 & 2 & 2 & 2 & 2 & 2 & 2 & 2 & 2 & 2 & 2 & 3 & 3 & 3 & 3 & 3 & 3 & 3 & 3 & 3 & 4 & 4 & 5 & 6 \\
   &   &   &   & $m_{53}$ & 2 & 2 & 2 & 2 & 2 & 3 & 3 & 3 & 3 & 3 & 4 & 4 & 4 & 4 & 4 & 5 & 5 & 5 & 5 & 5 & 6 & 6 & 6 & 6 & 6 & 7 & 7 & 7 & 7 & 7 & 2 & 3 & 3 & 4 & 4 & 5 & 5 & 6 & 6 & 3 & 4 & 3 & 3 \\
 $m_{61}$ & $m_{62}$ & $m_{63}$ & $m_{64}$ & $m_{54}$ & 2 & 3 & 4 & 5 & 6 & 2 & 3 & 4 & 5 & 6 & 2 & 3 & 4 & 5 & 6 & 2 & 3 & 4 & 5 & 6 & 2 & 3 & 4 & 5 & 6 & 2 & 3 & 4 & 5 & 6 & 3 & 2 & 3 & 2 & 3 & 2 & 3 & 2 & 3 & 2 & 2 & 2 & 2 \\
 3 & 2 & 2 & 2 &   & . & . & . & . & . & c & . & . & . & . & . & . & . & . & . & . & . & . & . & . & . & . & . & . & . & . & . & . & . & . & . & . & . & . & . & . & . & . & . & . & . & . & . \\
 3 & 2 & 2 & 3 &   & . & . & . & . & . & c & . & . & . & . & c & . & . & . & . & c & . & . & . & . & c & . & . & . & . & c & . & . & . & . & . & c & . & . & . & . & . & . & . & . & . & . & . \\
 3 & 2 & 2 & 4 &   & . & . & . & . & . & c & . & . & . & . & c & . & . & . & . & c & . & . & . & . & c & . & . & . & . & c & . & . & . & . & . & c & n & c & n & c & n & n & n & c & . & . & . \\
 3 & 2 & 2 & 5 &   & . & . & . & . & . & c & . & . & . & . & c & . & . & . & . & c & . & . & . & . & c & . & . & . & . & c & . & . & . & . & . & c & n & c & n & c & n & n & n & c & n & c & . \\
\end{tabular}
\par}
\caption{Cubes $(2,2,2,3,\ \ 2,\ast,\ast,\ast,\ \ 3,\ast,\ast,\ast)$}
\end{table}

\begin{table}[h!]
{\centering\setlength{\tabcolsep}{0.44pt}\footnotesize
\begin{tabular}{cccccccccccccccccccccccccccccccccccccccccccccccc}
 $m_{31}$ & $m_{32}$ & $m_{41}$ & $m_{42}$ & $m_{51}$ & 2 & 2 & 2 & 2 & 2 & 2 & 2 & 2 & 2 & 2 & 2 & 2 & 2 & 2 & 2 & 2 & 2 & 2 & 2 & 2 & 2 & 2 & 2 & 2 & 2 & 2 & 2 & 2 & 2 & 2 & 2 & 2 & 2 & 2 & 2 & 2 & 2 & 2 & 2 & 2 & 2 & 2 & 2 \\
 2 & 2 & 2 & 3 & $m_{52}$ & 2 & 2 & 2 & 2 & 2 & 2 & 2 & 2 & 2 & 2 & 2 & 2 & 2 & 2 & 2 & 2 & 2 & 2 & 2 & 2 & 2 & 2 & 2 & 2 & 2 & 2 & 2 & 2 & 2 & 2 & 3 & 3 & 3 & 3 & 3 & 3 & 3 & 3 & 3 & 4 & 4 & 5 & 6 \\
   &   &   &   & $m_{53}$ & 2 & 2 & 2 & 2 & 2 & 3 & 3 & 3 & 3 & 3 & 4 & 4 & 4 & 4 & 4 & 5 & 5 & 5 & 5 & 5 & 6 & 6 & 6 & 6 & 6 & 7 & 7 & 7 & 7 & 7 & 2 & 3 & 3 & 4 & 4 & 5 & 5 & 6 & 6 & 3 & 4 & 3 & 3 \\
 $m_{61}$ & $m_{62}$ & $m_{63}$ & $m_{64}$ & $m_{54}$ & 2 & 3 & 4 & 5 & 6 & 2 & 3 & 4 & 5 & 6 & 2 & 3 & 4 & 5 & 6 & 2 & 3 & 4 & 5 & 6 & 2 & 3 & 4 & 5 & 6 & 2 & 3 & 4 & 5 & 6 & 3 & 2 & 3 & 2 & 3 & 2 & 3 & 2 & 3 & 2 & 2 & 2 & 2 \\
 4 & 2 & 2 & 2 &   & . & c & c & c & n & c & c & c & c & n & c & . & . & . & . & . & . & . & . & . & . & . & . & . & . & . & . & . & . & . & . & . & . & . & . & . & . & . & . & . & . & . & . \\
 4 & 2 & 2 & 3 &   & . & c & c & c & n & c & c & c & c & n & c & c & c & c & n & c & c & c & c & n & c & c & c & c & n & c & c & c & c & n & . & c & n & c & . & . & . & . & . & . & . & . & . \\
 4 & 2 & 2 & 4 &   & . & n & n & n & n & n & n & n & n & n & n & n & n & n & n & n & n & n & n & n & n & n & n & n & n & n & n & n & n & n & . & n & n & n & n & n & n & n & n & n & n & . & . \\
 4 & 2 & 3 & 2 &   & . & c & c & c & n & c & c & c & c & n & c & c & c & c & n & c & c & c & c & n & c & c & c & c & n & c & c & c & c & n & . & c & n & c & n & c & n & n & n & c & n & c & n \\
\end{tabular}
\par}
\caption{Cubes $(2,2,2,3,\ \ 2,\ast,\ast,\ast,\ \ 4,\ast,\ast,\ast)$}
\end{table}

\begin{table}[h!]
{\centering\setlength{\tabcolsep}{0.44pt}\footnotesize
\begin{tabular}{cccccccccccccccccccccccccccccccccccccccccccccccc}
 $m_{31}$ & $m_{32}$ & $m_{41}$ & $m_{42}$ & $m_{51}$ & 2 & 2 & 2 & 2 & 2 & 2 & 2 & 2 & 2 & 2 & 2 & 2 & 2 & 2 & 2 & 2 & 2 & 2 & 2 & 2 & 2 & 2 & 2 & 2 & 2 & 2 & 2 & 2 & 2 & 2 & 2 & 2 & 2 & 2 & 2 & 2 & 2 & 2 & 2 & 2 & 2 & 2 & 2 \\
 2 & 2 & 2 & 3 & $m_{52}$ & 2 & 2 & 2 & 2 & 2 & 2 & 2 & 2 & 2 & 2 & 2 & 2 & 2 & 2 & 2 & 2 & 2 & 2 & 2 & 2 & 2 & 2 & 2 & 2 & 2 & 2 & 2 & 2 & 2 & 2 & 3 & 3 & 3 & 3 & 3 & 3 & 3 & 3 & 3 & 4 & 4 & 5 & 6 \\
   &   &   &   & $m_{53}$ & 2 & 2 & 2 & 2 & 2 & 3 & 3 & 3 & 3 & 3 & 4 & 4 & 4 & 4 & 4 & 5 & 5 & 5 & 5 & 5 & 6 & 6 & 6 & 6 & 6 & 7 & 7 & 7 & 7 & 7 & 2 & 3 & 3 & 4 & 4 & 5 & 5 & 6 & 6 & 3 & 4 & 3 & 3 \\
 $m_{61}$ & $m_{62}$ & $m_{63}$ & $m_{64}$ & $m_{54}$ & 2 & 3 & 4 & 5 & 6 & 2 & 3 & 4 & 5 & 6 & 2 & 3 & 4 & 5 & 6 & 2 & 3 & 4 & 5 & 6 & 2 & 3 & 4 & 5 & 6 & 2 & 3 & 4 & 5 & 6 & 3 & 2 & 3 & 2 & 3 & 2 & 3 & 2 & 3 & 2 & 2 & 2 & 2 \\
 5 & 2 & 2 & 2 &   & . & c & c & c & n & c & c & c & c & n & c & c & c & c & n & c & . & . & . & . & . & . & . & . & . & . & . & . & . & . & . & . & . & . & . & . & . & . & . & . & . & . & . \\
 5 & 2 & 2 & 3 &   & . & c & c & c & n & c & c & c & c & n & c & c & c & c & n & c & c & c & c & n & c & c & c & c & n & c & c & c & c & n & . & c & n & c & n & c & . & . & . & . & . & . & . \\
 5 & 2 & 3 & 2 &   & . & c & c & c & n & c & c & c & c & n & c & c & c & c & n & c & c & c & c & n & c & c & c & c & n & c & c & c & c & n & . & c & n & c & n & c & n & n & n & c & n & c & n \\
 5 & 2 & 3 & 3 &   & . & c & c & c & n & c & c & c & c & n & c & c & c & c & n & c & c & c & c & n & c & c & c & c & n & c & c & c & c & n & . & c & n & c & n & c & n & n & n & c & n & c & n \\
\end{tabular}
\par}
\caption{Cubes $(2,2,2,3,\ \ 2,\ast,\ast,\ast,\ \ 5+,\ast,\ast,\ast)$}
\end{table}

\begin{table}[h!]
{\centering\setlength{\tabcolsep}{0.44pt}\footnotesize
\begin{tabular}{ccccccccccccccccccccccccccccccccc}
 $m_{31}$ & $m_{32}$ & $m_{41}$ & $m_{42}$ & $m_{51}$ & 3 & 3 & 3 & 3 & 3 & 3 & 3 & 3 & 3 & 3 & 3 & 3 & 3 & 3 & 3 & 3 & 3 & 3 & 3 & 3 & 3 & 3 & 3 & 3 & 3 & 3 & 3 & 3 \\
 2 & 2 & 2 & 3 & $m_{52}$ & 2 & 2 & 2 & 2 & 2 & 2 & 2 & 2 & 2 & 2 & 2 & 2 & 2 & 2 & 2 & 2 & 2 & 2 & 2 & 2 & 3 & 3 & 3 & 3 & 3 & 3 & 3 & 4 \\
   &   &   &   & $m_{53}$ & 3 & 3 & 3 & 3 & 3 & 4 & 4 & 4 & 4 & 4 & 5 & 5 & 5 & 5 & 5 & 6 & 6 & 6 & 6 & 6 & 3 & 4 & 4 & 5 & 5 & 6 & 6 & 4 \\
 $m_{61}$ & $m_{62}$ & $m_{63}$ & $m_{64}$ & $m_{54}$ & 2 & 3 & 4 & 5 & 6 & 2 & 3 & 4 & 5 & 6 & 2 & 3 & 4 & 5 & 6 & 2 & 3 & 4 & 5 & 6 & 3 & 2 & 3 & 2 & 3 & 2 & 3 & 2 \\
 3 & 2 & 3 & 2 &   & c & . & . & . & . & . & . & . & . & . & . & . & . & . & . & . & . & . & . & . & . & . & . & . & . & . & . & . \\
 3 & 2 & 3 & 3 &   & c & c & . & . & . & . & . & . & . & . & . & . & . & . & . & . & . & . & . & . & . & . & . & . & . & . & . & . \\
 3 & 2 & 3 & 4 &   & c & c & c & . & . & . & . & . & . & . & . & . & . & . & . & . & . & . & . & . & . & . & . & . & . & . & . & . \\
 3 & 2 & 3 & 5 &   & c & c & c & c & . & . & . & . & . & . & . & . & . & . & . & . & . & . & . & . & . & . & . & . & . & . & . & . \\
\end{tabular}
\par}
\caption{Cubes $(2,2,2,3,\ \ 3,\ast,\ast,\ast,\ \ 3,\ast,\ast,\ast)$}
\end{table}

\begin{table}[h!]
{\centering\setlength{\tabcolsep}{0.44pt}\footnotesize
\begin{tabular}{ccccccccccccccccccccccccccccccccccccc}
 $m_{31}$ & $m_{32}$ & $m_{41}$ & $m_{42}$ & $m_{51}$ & 3 & 3 & 3 & 3 & 3 & 3 & 3 & 3 & 3 & 3 & 3 & 3 & 3 & 3 & 3 & 3 & 3 & 3 & 3 & 3 & 3 & 3 & 3 & 3 & 3 & 3 & 3 & 3 & 4 & 4 & 4 & 4 \\
 2 & 2 & 2 & 3 & $m_{52}$ & 2 & 2 & 2 & 2 & 2 & 2 & 2 & 2 & 2 & 2 & 2 & 2 & 2 & 2 & 2 & 2 & 2 & 2 & 2 & 2 & 3 & 3 & 3 & 3 & 3 & 3 & 3 & 4 & 2 & 2 & 2 & 3 \\
   &   &   &   & $m_{53}$ & 3 & 3 & 3 & 3 & 3 & 4 & 4 & 4 & 4 & 4 & 5 & 5 & 5 & 5 & 5 & 6 & 6 & 6 & 6 & 6 & 3 & 4 & 4 & 5 & 5 & 6 & 6 & 4 & 4 & 4 & 4 & 4 \\
 $m_{61}$ & $m_{62}$ & $m_{63}$ & $m_{64}$ & $m_{54}$ & 2 & 3 & 4 & 5 & 6 & 2 & 3 & 4 & 5 & 6 & 2 & 3 & 4 & 5 & 6 & 2 & 3 & 4 & 5 & 6 & 3 & 2 & 3 & 2 & 3 & 2 & 3 & 2 & 2 & 3 & 4 & 3 \\
 4 & 2 & 3 & 2 &   & . & c & c & c & n & c & . & . & . & . & . & . & . & . & . & . & . & . & . & . & . & . & . & . & . & . & . & . & . & . & . & . \\
 4 & 2 & 3 & 3 &   & . & c & c & c & n & c & c & c & c & n & c & c & c & c & n & n & n & n & n & n & . & c & . & . & . & . & . & . & . & . & . & . \\
 4 & 2 & 3 & 4 &   & . & n & n & n & n & n & n & n & n & n & n & n & n & n & n & n & n & n & n & n & . & n & n & n & n & n & n & n & . & . & . & . \\
 4 & 2 & 4 & 2 &   & n & n & n & n & n & n & n & n & n & n & n & n & n & n & n & n & n & n & n & n & n & n & n & n & n & n & n & n & n & . & . & . \\
\end{tabular}
\par}
\caption{Cubes $(2,2,2,3,\ \ 3+,\ast,\ast,\ast,\ \ 4+,\ast,\ast,\ast)$}
\end{table}

\begin{table}[h!]
{\centering\setlength{\tabcolsep}{0.44pt}\footnotesize
\begin{tabular}{cccccccccccccccccccccccccccccccccccccccccccccccc}
 $m_{31}$ & $m_{32}$ & $m_{41}$ & $m_{42}$ & $m_{51}$ & 2 & 2 & 2 & 2 & 2 & 2 & 2 & 2 & 2 & 2 & 2 & 2 & 2 & 2 & 2 & 2 & 2 & 2 & 2 & 2 & 2 & 2 & 2 & 2 & 3 & 3 & 3 & 3 & 3 & 3 & 3 & 3 & 3 & 3 & 3 & 3 & 3 & 3 & 3 & 3 & 4 & 4 & 4 \\
 2 & 2 & 2 & 4 & $m_{52}$ & 2 & 2 & 2 & 2 & 2 & 2 & 2 & 2 & 2 & 2 & 2 & 2 & 2 & 2 & 2 & 2 & 2 & 2 & 3 & 3 & 3 & 3 & 4 & 4 & 2 & 2 & 2 & 2 & 2 & 2 & 2 & 2 & 2 & 2 & 2 & 2 & 3 & 3 & 3 & 4 & 2 & 2 & 2 \\
   &   &   &   & $m_{53}$ & 2 & 2 & 2 & 3 & 3 & 3 & 4 & 4 & 4 & 5 & 5 & 5 & 6 & 6 & 6 & 7 & 7 & 7 & 3 & 4 & 5 & 6 & 3 & 4 & 3 & 3 & 3 & 4 & 4 & 4 & 5 & 5 & 5 & 6 & 6 & 6 & 4 & 5 & 6 & 4 & 4 & 4 & 4 \\
 $m_{61}$ & $m_{62}$ & $m_{63}$ & $m_{64}$ & $m_{54}$ & 2 & 3 & 4 & 2 & 3 & 4 & 2 & 3 & 4 & 2 & 3 & 4 & 2 & 3 & 4 & 2 & 3 & 4 & 2 & 2 & 2 & 2 & 2 & 2 & 2 & 3 & 4 & 2 & 3 & 4 & 2 & 3 & 4 & 2 & 3 & 4 & 2 & 2 & 2 & 2 & 2 & 3 & 4 \\
 2 & 3 & 3 & 2 &   & . & . & . & . & . & . & . & . & . & . & . & . & . & . & . & . & . & . & c & . & . & . & . & . & . & . & . & . & . & . & . & . & . & . & . & . & . & . & . & . & . & . & . \\
 2 & 3 & 4 & 2 &   & . & . & . & . & . & . & . & . & . & . & . & . & . & . & . & . & . & . & c & c & . & . & . & . & . & . & . & . & . & . & . & . & . & . & . & . & . & . & . & . & . & . & . \\
 2 & 3 & 5 & 2 &   & . & . & . & . & . & . & . & . & . & . & . & . & . & . & . & . & . & . & c & c & c & . & . & . & . & . & . & . & . & . & . & . & . & . & . & . & . & . & . & . & . & . & . \\
 2 & 3 & 6 & 2 &   & . & . & . & . & . & . & . & . & . & . & . & . & . & . & . & . & . & . & n & n & n & n & . & . & . & . & . & . & . & . & . & . & . & . & . & . & . & . & . & . & . & . & . \\
\end{tabular}
\par}
\caption{Cubes $(2,2,2,4,\ \ \ast,\ast,\ast,\ast,\ \ 3-,\ast,\ast,\ast)$}
\end{table}

\begin{table}[h!]
{\centering\setlength{\tabcolsep}{0.44pt}\footnotesize
\begin{tabular}{cccccccccccccccccccccccccccccccccccccccccccccccc}
 $m_{31}$ & $m_{32}$ & $m_{41}$ & $m_{42}$ & $m_{51}$ & 2 & 2 & 2 & 2 & 2 & 2 & 2 & 2 & 2 & 2 & 2 & 2 & 2 & 2 & 2 & 2 & 2 & 2 & 2 & 2 & 2 & 2 & 2 & 2 & 3 & 3 & 3 & 3 & 3 & 3 & 3 & 3 & 3 & 3 & 3 & 3 & 3 & 3 & 3 & 3 & 4 & 4 & 4 \\
 2 & 2 & 2 & 4 & $m_{52}$ & 2 & 2 & 2 & 2 & 2 & 2 & 2 & 2 & 2 & 2 & 2 & 2 & 2 & 2 & 2 & 2 & 2 & 2 & 3 & 3 & 3 & 3 & 4 & 4 & 2 & 2 & 2 & 2 & 2 & 2 & 2 & 2 & 2 & 2 & 2 & 2 & 3 & 3 & 3 & 4 & 2 & 2 & 2 \\
   &   &   &   & $m_{53}$ & 2 & 2 & 2 & 3 & 3 & 3 & 4 & 4 & 4 & 5 & 5 & 5 & 6 & 6 & 6 & 7 & 7 & 7 & 3 & 4 & 5 & 6 & 3 & 4 & 3 & 3 & 3 & 4 & 4 & 4 & 5 & 5 & 5 & 6 & 6 & 6 & 4 & 5 & 6 & 4 & 4 & 4 & 4 \\
 $m_{61}$ & $m_{62}$ & $m_{63}$ & $m_{64}$ & $m_{54}$ & 2 & 3 & 4 & 2 & 3 & 4 & 2 & 3 & 4 & 2 & 3 & 4 & 2 & 3 & 4 & 2 & 3 & 4 & 2 & 2 & 2 & 2 & 2 & 2 & 2 & 3 & 4 & 2 & 3 & 4 & 2 & 3 & 4 & 2 & 3 & 4 & 2 & 2 & 2 & 2 & 2 & 3 & 4 \\
 4 & 2 & 2 & 2 &   & . & . & . & . & . & . & c & . & . & . & . & . & . & . & . & . & . & . & . & . & . & . & . & . & . & . & . & . & . & . & . & . & . & . & . & . & . & . & . & . & . & . & . \\
 4 & 2 & 2 & 3 &   & . & . & . & c & . & . & c & . & . & c & . & . & c & . & . & c & . & . & c & c & . & . & . & . & . & . & . & . & . & . & . & . & . & . & . & . & . & . & . & . & . & . & . \\
 4 & 2 & 2 & 4 &   & . & . & . & n & . & . & n & . & . & n & . & . & n & . & . & n & . & . & n & n & n & n & n & n & . & . & . & . & . & . & . & . & . & . & . & . & . & . & . & . & . & . & . \\
 4 & 2 & 3 & 2 &   & . & c & . & c & c & . & c & c & . & . & c & . & . & c & . & . & c & . & c & c & c & n & n & n & . & c & n & c & . & . & . & . & . & . & . & . & . & . & . & . & . & . & . \\
\end{tabular}
\par}
\caption{Cubes $(2,2,2,4,\ \ \ast,\ast,\ast,\ast,\ \ 4+,\ast,\ast,\ast)$}
\end{table}

\begin{table}[h!]
{\centering\setlength{\tabcolsep}{0.44pt}\footnotesize
\begin{tabular}{cccccccccccccccccccccccccccccccccc}
 $m_{31}$ & $m_{32}$ & $m_{41}$ & $m_{42}$ & $m_{51}$ & 2 & 2 & 2 & 2 & 2 & 2 & 2 & 2 & 2 & 2 & 2 & 2 & 2 & 2 & 2 & 2 & 3 & 3 & 3 & 3 & 3 & 3 & 3 & 3 & 3 & 3 & 3 & 4 & 4 \\
 2 & 2 & 2 & 5 & $m_{52}$ & 2 & 2 & 2 & 2 & 2 & 2 & 2 & 2 & 2 & 2 & 2 & 2 & 3 & 3 & 3 & 3 & 2 & 2 & 2 & 2 & 2 & 2 & 2 & 2 & 3 & 3 & 3 & 2 & 2 \\
   &   &   &   & $m_{53}$ & 2 & 2 & 3 & 3 & 4 & 4 & 5 & 5 & 6 & 6 & 7 & 7 & 3 & 4 & 5 & 6 & 3 & 3 & 4 & 4 & 5 & 5 & 6 & 6 & 4 & 5 & 6 & 4 & 4 \\
 $m_{61}$ & $m_{62}$ & $m_{63}$ & $m_{64}$ & $m_{54}$ & 2 & 3 & 2 & 3 & 2 & 3 & 2 & 3 & 2 & 3 & 2 & 3 & 2 & 2 & 2 & 2 & 2 & 3 & 2 & 3 & 2 & 3 & 2 & 3 & 2 & 2 & 2 & 2 & 3 \\
 2 & 3 & 3 & 2 &   & . & . & . & . & . & . & . & . & . & . & . & . & c & . & . & . & . & . & . & . & . & . & . & . & . & . & . & . & . \\
 2 & 3 & 4 & 2 &   & . & . & . & . & . & . & . & . & . & . & . & . & c & c & . & . & . & . & . & . & . & . & . & . & . & . & . & . & . \\
 2 & 3 & 5 & 2 &   & . & . & . & . & . & . & . & . & . & . & . & . & c & c & c & . & . & . & . & . & . & . & . & . & . & . & . & . & . \\
 2 & 3 & 6 & 2 &   & . & . & . & . & . & . & . & . & . & . & . & . & n & n & n & n & . & . & . & . & . & . & . & . & . & . & . & . & . \\
\end{tabular}
\par}
\caption{Cubes $(2,2,2,5,\ \ \ast,\ast,\ast,\ast,\ \ \ast,\ast,\ast,\ast)$}
\end{table}

\begin{table}[h!]
{\centering\setlength{\tabcolsep}{0.44pt}\footnotesize
\begin{tabular}{cccccccccccccccccccccccccccccccccc}
 $m_{31}$ & $m_{32}$ & $m_{41}$ & $m_{42}$ & $m_{51}$ & 2 & 2 & 2 & 2 & 2 & 2 & 2 & 2 & 2 & 2 & 2 & 2 & 2 & 2 & 2 & 2 & 3 & 3 & 3 & 3 & 3 & 3 & 3 & 3 & 3 & 3 & 3 & 4 & 4 \\
 2 & 2 & 2 & 6 & $m_{52}$ & 2 & 2 & 2 & 2 & 2 & 2 & 2 & 2 & 2 & 2 & 2 & 2 & 3 & 3 & 3 & 3 & 2 & 2 & 2 & 2 & 2 & 2 & 2 & 2 & 3 & 3 & 3 & 2 & 2 \\
   &   &   &   & $m_{53}$ & 2 & 2 & 3 & 3 & 4 & 4 & 5 & 5 & 6 & 6 & 7 & 7 & 3 & 4 & 5 & 6 & 3 & 3 & 4 & 4 & 5 & 5 & 6 & 6 & 4 & 5 & 6 & 4 & 4 \\
 $m_{61}$ & $m_{62}$ & $m_{63}$ & $m_{64}$ & $m_{54}$ & 2 & 3 & 2 & 3 & 2 & 3 & 2 & 3 & 2 & 3 & 2 & 3 & 2 & 2 & 2 & 2 & 2 & 3 & 2 & 3 & 2 & 3 & 2 & 3 & 2 & 2 & 2 & 2 & 3 \\
 2 & 3 & 3 & 2 &   & . & . & . & . & . & . & . & . & . & . & . & . & n & . & . & . & . & . & . & . & . & . & . & . & . & . & . & . & . \\
 2 & 3 & 4 & 2 &   & . & . & . & . & . & . & . & . & . & . & . & . & n & n & . & . & . & . & . & . & . & . & . & . & . & . & . & . & . \\
 2 & 3 & 5 & 2 &   & . & . & . & . & . & . & . & . & . & . & . & . & n & n & n & . & . & . & . & . & . & . & . & . & . & . & . & . & . \\
 2 & 3 & 6 & 2 &   & . & . & . & . & . & . & . & . & . & . & . & . & n & n & n & n & . & . & . & . & . & . & . & . & . & . & . & . & . \\
\end{tabular}
\par}
\caption{Cubes $(2,2,2,6,\ \ \ast,\ast,\ast,\ast,\ \ \ast,\ast,\ast,\ast)$}
\end{table}

\begin{table}[h!]
{\centering\setlength{\tabcolsep}{0.44pt}\footnotesize
\begin{tabular}{ccccccccccc}
 $m_{31}$ & $m_{32}$ & $m_{41}$ & $m_{42}$ & $m_{51}$ & 2 & 3 & 3 & 3 & 3 & 4 \\
 2 & 2 & 2 & 7 & $m_{52}$ & 2 & 2 & 2 & 2 & 2 & 2 \\
   &   &   &   & $m_{53}$ & 7 & 3 & 4 & 5 & 6 & 4 \\
 $m_{61}$ & $m_{62}$ & $m_{63}$ & $m_{64}$ & $m_{54}$ & 2 & 2 & 2 & 2 & 2 & 2 \\
 3 & 2 & 3 & 2 &   & . & c & . & . & . & . \\
 3 & 2 & 4 & 2 &   & . & c & c & . & . & . \\
 3 & 2 & 5 & 2 &   & . & c & c & c & . & . \\
 3 & 2 & 6 & 2 &   & . & n & n & n & n & . \\
 4 & 2 & 3 & 2 &   & . & . & c & . & . & . \\
 4 & 2 & 4 & 2 &   & . & n & n & n & n & n \\
 5 & 2 & 3 & 2 &   & . & . & c & c & . & . \\
 6 & 2 & 3 & 2 &   & . & . & n & n & n & . \\
 7 & 2 & 2 & 2 &   & c & . & . & . & . & . \\
\end{tabular}
\par}
\caption{Cubes $(2,2,2,7,\ \ \ast,\ast,\ast,\ast,\ \ \ast,\ast,\ast,\ast)$}
\end{table}

\begin{table}[h!]
{\centering\setlength{\tabcolsep}{0.44pt}\footnotesize
\begin{tabular}{cccccccccccccccccccccccccccccccccccccccccccccccccccc}
 $m_{31}$ & $m_{32}$ & $m_{41}$ & $m_{42}$ & $m_{51}$ & 2 & 2 & 2 & 2 & 2 & 2 & 2 & 2 & 2 & 2 & 2 & 2 & 2 & 2 & 2 & 2 & 2 & 2 & 2 & 2 & 2 & 2 & 2 & 2 & 2 & 2 & 2 & 2 & 2 & 2 & 2 & 2 & 2 & 2 & 2 & 2 & 2 & 2 & 2 & 2 & 2 & 2 & 2 & 2 & 2 & 2 & 2 \\
 2 & 2 & 3 & 3 & $m_{52}$ & 2 & 2 & 2 & 2 & 2 & 2 & 2 & 2 & 2 & 2 & 2 & 2 & 2 & 2 & 2 & 2 & 2 & 2 & 2 & 2 & 2 & 2 & 2 & 2 & 2 & 2 & 2 & 2 & 2 & 2 & 3 & 3 & 3 & 3 & 3 & 3 & 3 & 3 & 3 & 3 & 4 & 4 & 4 & 5 & 5 & 6 & 6 \\
   &   &   &   & $m_{53}$ & 2 & 2 & 2 & 2 & 2 & 3 & 3 & 3 & 3 & 3 & 4 & 4 & 4 & 4 & 4 & 5 & 5 & 5 & 5 & 5 & 6 & 6 & 6 & 6 & 6 & 7 & 7 & 7 & 7 & 7 & 2 & 2 & 3 & 3 & 4 & 4 & 5 & 5 & 6 & 6 & 2 & 3 & 4 & 2 & 3 & 2 & 3 \\
 $m_{61}$ & $m_{62}$ & $m_{63}$ & $m_{64}$ & $m_{54}$ & 2 & 3 & 4 & 5 & 6 & 2 & 3 & 4 & 5 & 6 & 2 & 3 & 4 & 5 & 6 & 2 & 3 & 4 & 5 & 6 & 2 & 3 & 4 & 5 & 6 & 2 & 3 & 4 & 5 & 6 & 2 & 3 & 2 & 3 & 2 & 3 & 2 & 3 & 2 & 3 & 2 & 2 & 2 & 2 & 2 & 2 & 2 \\
 2 & 3 & 3 & 3 &   & . & . & . & . & . & . & . & . & . & . & . & . & . & . & . & . & . & . & . & . & . & . & . & . & . & . & . & . & . & . & n & . & n & n & . & . & . & . & . & . & . & . & . & . & . & . & . \\
 2 & 3 & 4 & 2 &   & . & . & . & . & . & . & . & . & . & . & . & . & . & . & . & . & . & . & . & . & . & . & . & . & . & . & . & . & . & . & . & n & . & n & . & . & . & . & . & . & . & . & . & . & . & . & . \\
 2 & 3 & 4 & 3 &   & . & . & . & . & . & . & . & . & . & . & . & . & . & . & . & . & . & . & . & . & . & . & . & . & . & . & . & . & . & . & n & n & n & n & . & n & . & . & . & . & . & . & . & . & . & . & . \\
 2 & 3 & 5 & 2 &   & . & . & . & . & . & . & . & . & . & . & . & . & . & . & . & . & . & . & . & . & . & . & . & . & . & . & . & . & . & . & . & n & . & n & . & . & . & . & . & . & . & . & . & . & . & . & . \\
\end{tabular}
\par}
\caption{Cubes $(2,2,3,3,\ \ 2,\ast,\ast,\ast,\ \ 3-,\ast,\ast,\ast)$}
\end{table}

\begin{table}[h!]
{\centering\setlength{\tabcolsep}{0.44pt}\footnotesize
\begin{tabular}{cccccccccccccccccccccccccccccccccccccccccccccccccccc}
 $m_{31}$ & $m_{32}$ & $m_{41}$ & $m_{42}$ & $m_{51}$ & 2 & 2 & 2 & 2 & 2 & 2 & 2 & 2 & 2 & 2 & 2 & 2 & 2 & 2 & 2 & 2 & 2 & 2 & 2 & 2 & 2 & 2 & 2 & 2 & 2 & 2 & 2 & 2 & 2 & 2 & 2 & 2 & 2 & 2 & 2 & 2 & 2 & 2 & 2 & 2 & 2 & 2 & 2 & 2 & 2 & 2 & 2 \\
 2 & 2 & 3 & 3 & $m_{52}$ & 2 & 2 & 2 & 2 & 2 & 2 & 2 & 2 & 2 & 2 & 2 & 2 & 2 & 2 & 2 & 2 & 2 & 2 & 2 & 2 & 2 & 2 & 2 & 2 & 2 & 2 & 2 & 2 & 2 & 2 & 3 & 3 & 3 & 3 & 3 & 3 & 3 & 3 & 3 & 3 & 4 & 4 & 4 & 5 & 5 & 6 & 6 \\
   &   &   &   & $m_{53}$ & 2 & 2 & 2 & 2 & 2 & 3 & 3 & 3 & 3 & 3 & 4 & 4 & 4 & 4 & 4 & 5 & 5 & 5 & 5 & 5 & 6 & 6 & 6 & 6 & 6 & 7 & 7 & 7 & 7 & 7 & 2 & 2 & 3 & 3 & 4 & 4 & 5 & 5 & 6 & 6 & 2 & 3 & 4 & 2 & 3 & 2 & 3 \\
 $m_{61}$ & $m_{62}$ & $m_{63}$ & $m_{64}$ & $m_{54}$ & 2 & 3 & 4 & 5 & 6 & 2 & 3 & 4 & 5 & 6 & 2 & 3 & 4 & 5 & 6 & 2 & 3 & 4 & 5 & 6 & 2 & 3 & 4 & 5 & 6 & 2 & 3 & 4 & 5 & 6 & 2 & 3 & 2 & 3 & 2 & 3 & 2 & 3 & 2 & 3 & 2 & 2 & 2 & 2 & 2 & 2 & 2 \\
 4 & 2 & 2 & 2 &   & . & . & . & . & . & . & . & . & . & . & . & . & . & . & . & . & . & . & . & . & . & . & . & . & . & . & . & . & . & . & . & . & . & . & . & n & . & n & . & n & . & . & . & . & . & . & . \\
 4 & 2 & 3 & 2 &   & . & . & . & . & . & . & . & . & . & . & . & . & . & . & . & . & . & . & . & . & . & . & . & . & . & . & . & . & . & . & . & n & c & n & c & n & c & n & n & n & . & c & . & . & . & . & . \\
 4 & 2 & 4 & 2 &   & . & . & . & . & . & . & . & . & . & . & . & . & . & . & . & . & . & . & . & . & . & . & . & . & . & . & . & . & . & . & . & n & n & n & n & n & n & n & n & n & . & n & n & . & . & . & . \\
 4 & 3 & 2 & 2 &   & . & . & . & . & . & . & . & . & . & . & . & . & . & . & . & . & . & . & . & . & . & . & . & . & . & . & . & . & . & . & . & . & . & . & . & n & . & n & . & n & . & . & . & . & . & . & . \\
\end{tabular}
\par}
\caption{Cubes $(2,2,3,3,\ \ 2,\ast,\ast,\ast,\ \ 4+,\ast,\ast,\ast)$}
\end{table}

\begin{table}[h!]
{\centering\setlength{\tabcolsep}{0.44pt}\footnotesize
\begin{tabular}{cccccccccccccccccccccccccc}
 $m_{31}$ & $m_{32}$ & $m_{41}$ & $m_{42}$ & $m_{51}$ & 3 & 3 & 3 & 3 & 3 & 3 & 3 & 3 & 3 & 3 & 3 & 3 & 3 & 3 & 4 & 4 & 4 & 4 & 5 & 5 & 6 \\
 2 & 2 & 3 & 3 & $m_{52}$ & 3 & 3 & 3 & 3 & 3 & 3 & 3 & 3 & 3 & 3 & 4 & 4 & 5 & 6 & 4 & 4 & 5 & 6 & 5 & 6 & 6 \\
   &   &   &   & $m_{53}$ & 2 & 2 & 3 & 3 & 4 & 4 & 5 & 5 & 6 & 6 & 3 & 4 & 3 & 3 & 3 & 4 & 3 & 3 & 3 & 3 & 3 \\
 $m_{61}$ & $m_{62}$ & $m_{63}$ & $m_{64}$ & $m_{54}$ & 2 & 3 & 2 & 3 & 2 & 3 & 2 & 3 & 2 & 3 & 2 & 2 & 2 & 2 & 2 & 2 & 2 & 2 & 2 & 2 & 2 \\
 3 & 3 & 2 & 3 &   & . & n & . & . & . & . & . & . & . & . & . & . & . & . & . & . & . & . & . & . & . \\
 3 & 3 & 3 & 2 &   & . & n & c & . & . & . & . & . & . & . & . & . & . & . & . & . & . & . & . & . & . \\
 3 & 3 & 3 & 3 &   & . & n & n & n & . & . & . & . & . & . & . & . & . & . & . & . & . & . & . & . & . \\
 3 & 3 & 4 & 2 &   & . & n & c & n & c & . & . & . & . & . & . & . & . & . & . & . & . & . & . & . & . \\
\end{tabular}
\par}
\caption{Cubes $(2,2,3,3,\ \ 3+,\ast,\ast,\ast,\ \ \ast,\ast,\ast,\ast)$}
\end{table}

\begin{table}[h!]
{\centering\setlength{\tabcolsep}{0.44pt}\footnotesize
\begin{tabular}{ccccccccccccccccccccccccccccccccccccccccccccccccc}
 $m_{31}$ & $m_{32}$ & $m_{41}$ & $m_{42}$ & $m_{51}$ & 2 & 2 & 2 & 2 & 2 & 2 & 2 & 2 & 2 & 2 & 2 & 2 & 2 & 2 & 2 & 2 & 2 & 2 & 2 & 2 & 2 & 2 & 2 & 2 & 2 & 3 & 3 & 3 & 3 & 3 & 3 & 3 & 3 & 3 & 3 & 3 & 3 & 3 & 3 & 3 & 4 & 4 & 4 & 4 \\
 2 & 2 & 3 & 4 & $m_{52}$ & 2 & 2 & 2 & 2 & 2 & 2 & 2 & 2 & 2 & 2 & 2 & 2 & 2 & 2 & 2 & 2 & 2 & 3 & 3 & 3 & 3 & 3 & 4 & 4 & 4 & 2 & 2 & 2 & 2 & 2 & 2 & 2 & 2 & 2 & 3 & 3 & 3 & 3 & 4 & 4 & 2 & 3 & 4 & 4 \\
   &   &   &   & $m_{53}$ & 2 & 2 & 3 & 3 & 3 & 4 & 4 & 4 & 5 & 5 & 5 & 6 & 6 & 6 & 7 & 7 & 7 & 2 & 3 & 4 & 5 & 6 & 2 & 3 & 4 & 2 & 3 & 3 & 4 & 4 & 5 & 5 & 6 & 6 & 3 & 4 & 5 & 6 & 3 & 4 & 4 & 4 & 3 & 4 \\
 $m_{61}$ & $m_{62}$ & $m_{63}$ & $m_{64}$ & $m_{54}$ & 3 & 4 & 2 & 3 & 4 & 2 & 3 & 4 & 2 & 3 & 4 & 2 & 3 & 4 & 2 & 3 & 4 & 2 & 2 & 2 & 2 & 2 & 2 & 2 & 2 & 3 & 2 & 3 & 2 & 3 & 2 & 3 & 2 & 3 & 2 & 2 & 2 & 2 & 2 & 2 & 2 & 2 & 2 & 2 \\
 2 & 4 & 4 & 2 &   & . & . & . & . & . & . & . & . & . & . & . & . & . & . & . & . & . & . & . & . & . & . & . & . & n & . & . & . & . & . & . & . & . & . & . & . & . & . & . & . & . & . & . & . \\
 3 & 2 & 2 & 3 &   & . & . & . & . & . & . & . & . & . & . & . & . & . & . & . & . & . & . & n & n & n & n & . & n & n & . & . & . & . & . & . & . & . & . & . & . & . & . & . & . & . & . & . & . \\
 3 & 2 & 3 & 2 &   & . & . & . & . & . & . & . & . & . & . & . & . & . & . & . & . & . & . & . & c & c & n & . & . & n & . & . & . & . & . & . & . & . & . & . & . & . & . & . & . & . & . & . & . \\
 3 & 2 & 3 & 3 &   & . & . & . & . & . & . & . & . & . & . & . & . & . & . & . & . & . & . & n & n & n & n & . & n & n & . & . & . & . & . & . & . & . & . & . & . & . & . & . & . & . & . & . & . \\
\end{tabular}
\par}
\caption{Cubes $(2,2,3,4,\ \ \ast,\ast,\ast,\ast,\ \ \ast,\ast,\ast,\ast)$}
\end{table}

\begin{table}[h!]
{\centering\setlength{\tabcolsep}{0.44pt}\footnotesize
\begin{tabular}{cccccccccccccccccccccccccccccccccc}
 $m_{31}$ & $m_{32}$ & $m_{41}$ & $m_{42}$ & $m_{51}$ & 2 & 2 & 2 & 2 & 2 & 2 & 2 & 2 & 2 & 2 & 2 & 2 & 2 & 2 & 3 & 3 & 3 & 3 & 3 & 3 & 3 & 3 & 3 & 3 & 3 & 3 & 3 & 4 & 4 \\
 2 & 2 & 3 & 5 & $m_{52}$ & 2 & 2 & 2 & 2 & 2 & 2 & 2 & 2 & 2 & 3 & 3 & 3 & 3 & 3 & 2 & 2 & 2 & 2 & 2 & 2 & 2 & 2 & 2 & 3 & 3 & 3 & 3 & 2 & 3 \\
   &   &   &   & $m_{53}$ & 2 & 3 & 4 & 5 & 5 & 6 & 6 & 7 & 7 & 2 & 3 & 4 & 5 & 6 & 2 & 3 & 3 & 4 & 4 & 5 & 5 & 6 & 6 & 3 & 4 & 5 & 6 & 4 & 4 \\
 $m_{61}$ & $m_{62}$ & $m_{63}$ & $m_{64}$ & $m_{54}$ & 3 & 3 & 3 & 2 & 3 & 2 & 3 & 2 & 3 & 2 & 2 & 2 & 2 & 2 & 3 & 2 & 3 & 2 & 3 & 2 & 3 & 2 & 3 & 2 & 2 & 2 & 2 & 2 & 2 \\
 3 & 2 & 2 & 3 &   & . & . & . & . & . & . & . & . & . & . & n & n & n & n & . & . & . & . & . & . & . & . & . & . & . & . & . & . & . \\
 3 & 2 & 3 & 2 &   & . & . & . & . & . & . & . & . & . & . & . & . & c & n & . & . & . & . & . & . & . & . & . & . & . & . & . & . & . \\
 3 & 2 & 3 & 3 &   & . & . & . & . & . & . & . & . & . & . & n & n & n & n & . & . & . & . & . & . & . & . & . & . & . & . & . & . & . \\
 3 & 2 & 4 & 2 &   & . & . & . & . & . & . & . & . & . & . & . & c & c & n & . & . & . & . & . & . & . & . & . & . & . & . & . & . & . \\
\end{tabular}
\par}
\caption{Cubes $(2,2,3,5,\ \ \ast,\ast,\ast,\ast,\ \ \ast,\ast,\ast,\ast)$}
\end{table}

\begin{table}[h!]
{\centering\setlength{\tabcolsep}{0.44pt}\footnotesize
\begin{tabular}{cccccccccccccccccccccccccccccccc}
 $m_{31}$ & $m_{32}$ & $m_{41}$ & $m_{42}$ & $m_{51}$ & 2 & 2 & 2 & 2 & 2 & 2 & 2 & 2 & 2 & 2 & 2 & 2 & 2 & 3 & 3 & 3 & 3 & 3 & 3 & 3 & 3 & 3 & 3 & 3 & 3 & 4 & 4 \\
 2 & 2 & 3 & 6 & $m_{52}$ & 2 & 2 & 2 & 2 & 2 & 2 & 2 & 2 & 3 & 3 & 3 & 3 & 3 & 2 & 2 & 2 & 2 & 2 & 2 & 2 & 2 & 3 & 3 & 3 & 3 & 2 & 3 \\
   &   &   &   & $m_{53}$ & 2 & 3 & 4 & 5 & 6 & 6 & 7 & 7 & 2 & 3 & 4 & 5 & 6 & 2 & 3 & 4 & 4 & 5 & 5 & 6 & 6 & 3 & 4 & 5 & 6 & 4 & 4 \\
 $m_{61}$ & $m_{62}$ & $m_{63}$ & $m_{64}$ & $m_{54}$ & 3 & 3 & 3 & 3 & 2 & 3 & 2 & 3 & 2 & 2 & 2 & 2 & 2 & 3 & 3 & 2 & 3 & 2 & 3 & 2 & 3 & 2 & 2 & 2 & 2 & 2 & 2 \\
 3 & 2 & 2 & 3 &   & . & . & . & . & . & . & . & . & . & n & n & n & n & . & . & . & . & . & . & . & . & . & . & . & . & . & . \\
 3 & 2 & 3 & 2 &   & . & . & . & . & . & . & . & . & . & . & . & . & n & . & . & . & . & . & . & . & . & . & . & . & . & . & . \\
 3 & 2 & 3 & 3 &   & . & . & . & . & . & . & . & . & . & n & n & n & n & . & . & . & . & . & . & . & . & . & . & . & . & . & . \\
 3 & 2 & 4 & 2 &   & . & . & . & . & . & . & . & . & . & . & n & n & n & . & . & . & . & . & . & . & . & . & . & . & . & . & . \\
\end{tabular}
\par}
\caption{Cubes $(2,2,3,6,\ \ \ast,\ast,\ast,\ast,\ \ \ast,\ast,\ast,\ast)$}
\end{table}

\begin{table}[h!]
{\centering\setlength{\tabcolsep}{0.44pt}\footnotesize
\begin{tabular}{cccccccccccccccccccccccccccccc}
 $m_{31}$ & $m_{32}$ & $m_{41}$ & $m_{42}$ & $m_{51}$ & 2 & 2 & 2 & 2 & 2 & 2 & 2 & 2 & 2 & 2 & 2 & 2 & 2 & 2 & 2 & 2 & 2 & 2 & 2 & 2 & 3 & 3 & 3 & 3 & 4 \\
 2 & 2 & 4 & 4 & $m_{52}$ & 2 & 2 & 2 & 2 & 2 & 2 & 2 & 2 & 2 & 2 & 2 & 2 & 2 & 2 & 2 & 2 & 3 & 3 & 3 & 4 & 3 & 3 & 3 & 4 & 4 \\
   &   &   &   & $m_{53}$ & 2 & 2 & 3 & 3 & 4 & 4 & 4 & 5 & 5 & 5 & 6 & 6 & 6 & 7 & 7 & 7 & 4 & 5 & 6 & 4 & 4 & 5 & 6 & 4 & 4 \\
 $m_{61}$ & $m_{62}$ & $m_{63}$ & $m_{64}$ & $m_{54}$ & 3 & 4 & 3 & 4 & 2 & 3 & 4 & 2 & 3 & 4 & 2 & 3 & 4 & 2 & 3 & 4 & 2 & 2 & 2 & 2 & 2 & 2 & 2 & 2 & 2 \\
 3 & 2 & 4 & 2 &   & . & . & . & . & . & . & . & . & . & . & . & . & . & . & . & . & c & . & . & . & . & . & . & . & . \\
 3 & 2 & 5 & 2 &   & . & . & . & . & . & . & . & . & . & . & . & . & . & . & . & . & c & c & . & . & . & . & . & . & . \\
 3 & 2 & 6 & 2 &   & . & . & . & . & . & . & . & . & . & . & . & . & . & . & . & . & n & n & n & . & . & . & . & . & . \\
 3 & 3 & 4 & 2 &   & . & . & . & . & . & . & . & . & . & . & . & . & . & . & . & . & c & c & n & . & c & . & . & . & . \\
 3 & 3 & 5 & 2 &   & . & . & . & . & . & . & . & . & . & . & . & . & . & . & . & . & c & c & n & n & c & c & . & . & . \\
 3 & 3 & 6 & 2 &   & . & . & . & . & . & . & . & . & . & . & . & . & . & . & . & . & n & n & n & n & n & n & n & . & . \\
 3 & 4 & 4 & 2 &   & . & . & . & . & . & . & . & . & . & . & . & . & . & . & . & . & n & n & n & . & n & n & n & . & . \\
 4 & 2 & 4 & 2 &   & . & . & . & . & . & . & . & . & . & . & . & . & . & . & . & . & n & n & n & n & . & . & . & . & . \\
 4 & 3 & 4 & 2 &   & . & . & . & . & . & . & . & . & . & . & . & . & . & . & . & . & n & n & n & n & . & . & . & n & . \\
 4 & 4 & 3 & 2 &   & n & . & n & . & . & n & . & . & n & . & . & n & . & . & n & . & . & . & . & . & . & . & . & . & . \\
 4 & 4 & 4 & 2 &   & n & n & n & n & n & n & n & n & n & n & n & n & n & n & n & n & n & n & n & n & n & n & n & n & n \\
\end{tabular}
\par}
\caption{Cubes $(2,2,4,4,\ \ \ast,\ast,\ast,\ast,\ \ \ast,\ast,\ast,\ast)$}
\end{table}

\begin{table}[h!]
{\centering\setlength{\tabcolsep}{0.44pt}\footnotesize
\begin{tabular}{cccccccccccccc}
 $m_{31}$ & $m_{32}$ & $m_{41}$ & $m_{42}$ & $m_{51}$ & 2 & 2 & 2 & 3 & 3 & 3 & 3 & 3 & 3 \\
 2 & 2 & 4 & 5 & $m_{52}$ & 3 & 3 & 3 & 2 & 2 & 2 & 3 & 3 & 3 \\
   &   &   &   & $m_{53}$ & 4 & 5 & 6 & 4 & 5 & 6 & 4 & 5 & 6 \\
 $m_{61}$ & $m_{62}$ & $m_{63}$ & $m_{64}$ & $m_{54}$ & 2 & 2 & 2 & 2 & 2 & 2 & 2 & 2 & 2 \\
 3 & 2 & 4 & 2 &   & . & c & n & . & . & . & . & . & . \\
 3 & 2 & 5 & 2 &   & c & c & n & . & . & . & . & . & . \\
 3 & 2 & 6 & 2 &   & n & n & n & . & . & . & . & . & . \\
 3 & 3 & 4 & 2 &   & . & c & n & . & c & n & . & . & . \\
 3 & 3 & 5 & 2 &   & c & c & n & . & c & n & c & c & . \\
 3 & 3 & 6 & 2 &   & n & n & n & n & n & n & n & n & n \\
 4 & 2 & 4 & 2 &   & . & n & n & . & . & . & . & . & n \\
 4 & 3 & 4 & 2 &   & . & n & n & . & n & n & . & n & n \\
\end{tabular}
\par}
\caption{Cubes $(2,2,4,5,\ \ \ast,\ast,\ast,\ast,\ \ \ast,\ast,\ast,\ast)$}
\end{table}

\begin{table}[h!]
{\centering\setlength{\tabcolsep}{0.44pt}\footnotesize
\begin{tabular}{ccccccccccccc}
 $m_{31}$ & $m_{32}$ & $m_{41}$ & $m_{42}$ & $m_{51}$ & 2 & 2 & 2 & 3 & 3 & 3 & 3 & 3 \\
 2 & 2 & 4 & 6 & $m_{52}$ & 3 & 3 & 3 & 2 & 2 & 3 & 3 & 3 \\
   &   &   &   & $m_{53}$ & 4 & 5 & 6 & 5 & 6 & 4 & 5 & 6 \\
 $m_{61}$ & $m_{62}$ & $m_{63}$ & $m_{64}$ & $m_{54}$ & 2 & 2 & 2 & 2 & 2 & 2 & 2 & 2 \\
 3 & 2 & 4 & 2 &   & . & . & n & . & . & . & . & . \\
 3 & 2 & 5 & 2 &   & . & n & n & . & . & . & . & . \\
 3 & 2 & 6 & 2 &   & n & n & n & . & . & . & . & . \\
 3 & 3 & 4 & 2 &   & . & . & n & . & n & . & . & . \\
 3 & 3 & 5 & 2 &   & . & n & n & n & n & . & n & . \\
 3 & 3 & 6 & 2 &   & n & n & n & n & n & n & n & n \\
 4 & 2 & 4 & 2 &   & . & . & n & . & . & . & . & . \\
 4 & 3 & 4 & 2 &   & . & . & n & . & n & . & . & n \\
\end{tabular}
\par}
\caption{Cubes $(2,2,4,6,\ \ \ast,\ast,\ast,\ast,\ \ \ast,\ast,\ast,\ast)$}
\end{table}

\begin{table}[h!]
{\centering\setlength{\tabcolsep}{0.44pt}\footnotesize
\begin{tabular}{cccccccccccccccc}
 $m_{31}$ & $m_{32}$ & $m_{41}$ & $m_{42}$ & $m_{51}$ & 2 & 2 & 2 & 2 & 2 & 2 & 2 & 2 & 2 & 2 & 3 \\
 2 & 3 & 3 & 2 & $m_{52}$ & 3 & 3 & 3 & 3 & 3 & 3 & 3 & 3 & 3 & 4 & 3 \\
   &   &   &   & $m_{53}$ & 2 & 2 & 2 & 2 & 3 & 3 & 3 & 3 & 3 & 2 & 3 \\
 $m_{61}$ & $m_{62}$ & $m_{63}$ & $m_{64}$ & $m_{54}$ & 3 & 4 & 5 & 6 & 2 & 3 & 4 & 5 & 6 & 4 & 3 \\
 3 & 2 & 2 & 3 &   & . & . & . & . & n & . & . & . & . & . & . \\
 3 & 2 & 3 & 2 &   & c & . & . & . & . & . & . & . & . & . & . \\
 3 & 2 & 3 & 3 &   & n & . & . & . & n & n & . & . & . & . & . \\
 3 & 2 & 4 & 2 &   & c & c & . & . & . & . & . & . & . & . & . \\
 3 & 2 & 4 & 3 &   & n & n & . & . & n & n & n & . & . & . & . \\
 3 & 2 & 5 & 2 &   & c & c & c & . & . & . & . & . & . & . & . \\
 3 & 2 & 5 & 3 &   & n & n & n & . & n & n & n & n & . & . & . \\
 3 & 2 & 6 & 2 &   & n & n & n & n & . & . & . & . & . & . & . \\
 3 & 2 & 6 & 3 &   & n & n & n & n & n & n & n & n & n & . & . \\
 3 & 3 & 2 & 3 &   & . & . & . & . & . & n & . & . & . & . & . \\
 3 & 3 & 3 & 2 &   & . & n & n & n & . & n & . & . & . & . & . \\
 3 & 3 & 3 & 3 &   & n & n & n & n & n & n & n & n & n & n & n \\
 3 & 4 & 2 & 3 &   & . & . & . & . & . & n & n & . & . & . & . \\
 3 & 5 & 2 & 3 &   & . & . & . & . & . & n & n & n & . & . & . \\
 3 & 6 & 2 & 3 &   & . & . & . & . & . & n & n & n & n & . & . \\
 4 & 2 & 3 & 2 &   & . & c & . & . & . & . & . & . & . & . & . \\
 4 & 2 & 4 & 2 &   & n & n & n & n & . & n & n & n & n & n & . \\
 4 & 3 & 3 & 2 &   & . & n & n & n & . & n & n & . & . & . & . \\
 5 & 2 & 3 & 2 &   & . & c & c & . & . & . & . & . & . & . & . \\
 5 & 3 & 3 & 2 &   & . & n & n & n & . & n & n & n & . & . & . \\
 6 & 2 & 3 & 2 &   & . & n & n & n & . & . & . & . & . & . & . \\
 6 & 3 & 3 & 2 &   & . & n & n & n & . & n & n & n & n & . & . \\
\end{tabular}
\par}
\caption{Cubes $(2,3,3,2,\ \ \ast,\ast,\ast,\ast,\ \ \ast,\ast,\ast,\ast)$}
\end{table}

\begin{table}[h!]
{\centering\setlength{\tabcolsep}{0.44pt}\footnotesize
\begin{tabular}{ccccccccccccc}
 $m_{31}$ & $m_{32}$ & $m_{41}$ & $m_{42}$ & $m_{51}$ & 2 & 2 & 2 & 3 & 3 & 3 & 3 & 3 \\
 2 & 3 & 3 & 3 & $m_{52}$ & 3 & 3 & 3 & 2 & 2 & 2 & 2 & 3 \\
   &   &   &   & $m_{53}$ & 2 & 3 & 3 & 3 & 4 & 5 & 6 & 3 \\
 $m_{61}$ & $m_{62}$ & $m_{63}$ & $m_{64}$ & $m_{54}$ & 3 & 2 & 3 & 3 & 3 & 3 & 3 & 3 \\
 3 & 3 & 3 & 2 &   & . & . & n & n & . & . & . & . \\
 3 & 3 & 3 & 3 &   & n & n & n & n & n & n & n & n \\
 4 & 2 & 3 & 2 &   & . & . & n & . & . & . & . & . \\
 4 & 2 & 4 & 2 &   & n & . & n & . & . & . & . & n \\
 4 & 3 & 3 & 2 &   & . & . & n & n & n & . & . & . \\
 5 & 2 & 3 & 2 &   & . & . & n & . & . & . & . & . \\
 5 & 3 & 3 & 2 &   & . & . & n & n & n & n & . & . \\
 6 & 2 & 3 & 2 &   & . & . & n & . & . & . & . & . \\
 6 & 3 & 3 & 2 &   & . & . & n & n & n & n & n & . \\
\end{tabular}
\par}
\caption{Cubes $(2,3,3,3,\ \ \ast,\ast,\ast,\ast,\ \ \ast,\ast,\ast,\ast)$}
\end{table}

\begin{table}[h!]
{\centering\setlength{\tabcolsep}{0.44pt}\footnotesize
\begin{tabular}{cccccccccccc}
 $m_{31}$ & $m_{32}$ & $m_{41}$ & $m_{42}$ & $m_{51}$ & 2 & 2 & 2 & 2 & 2 & 2 & 2 \\
 2 & 3 & 4 & 2 & $m_{52}$ & 3 & 3 & 3 & 4 & 4 & 5 & 6 \\
   &   &   &   & $m_{53}$ & 2 & 2 & 3 & 2 & 2 & 2 & 2 \\
 $m_{61}$ & $m_{62}$ & $m_{63}$ & $m_{64}$ & $m_{54}$ & 3 & 4 & 4 & 3 & 4 & 3 & 3 \\
 3 & 2 & 4 & 2 &   & . & . & . & c & . & . & . \\
 3 & 2 & 5 & 2 &   & . & . & . & c & n & c & . \\
 3 & 2 & 6 & 2 &   & . & . & . & n & n & n & n \\
 4 & 2 & 4 & 2 &   & n & n & n & n & n & . & . \\
 4 & 3 & 3 & 2 &   & . & n & n & . & . & . & . \\
\end{tabular}
\par}
\caption{Cubes $(2,3,4,2,\ \ \ast,\ast,\ast,\ast,\ \ \ast,\ast,\ast,\ast)$}
\end{table}

\begin{table}[h!]
{\centering\setlength{\tabcolsep}{0.44pt}\footnotesize
\begin{tabular}{cccccccccccccc}
 $m_{31}$ & $m_{32}$ & $m_{41}$ & $m_{42}$ & $m_{51}$ & $m_{52}$ & $m_{53}$ & $m_{54}$ & $m_{61}$ & $m_{62}$ & $m_{63}$ & $m_{64}$ &   & \text{type } \\
 2 & 2 & 3 & 7 & 4 & 2 & 4 & 2 & 4 & 2 & 4 & 2 &   & n \\
 2 & 2 & 5 & 5 & 2 & 3 & 5 & 2 & 3 & 2 & 5 & 2 &   & c \\
 2 & 2 & 5 & 5 & 2 & 3 & 5 & 2 & 3 & 2 & 6 & 2 &   & n \\
 2 & 2 & 5 & 5 & 2 & 3 & 6 & 2 & 3 & 2 & 6 & 2 &   & n \\
 2 & 2 & 5 & 5 & 2 & 3 & 5 & 2 & 3 & 3 & 5 & 2 &   & c \\
 2 & 2 & 5 & 5 & 2 & 3 & 5 & 2 & 3 & 3 & 6 & 2 &   & n \\
 2 & 2 & 5 & 5 & 2 & 3 & 6 & 2 & 3 & 3 & 5 & 2 &   & n \\
\end{tabular}
\par}
\caption{Remaining cubes}
\end{table}
\end{document}